\documentclass{amsart}
\usepackage[dvips]{epsfig}
\usepackage{graphicx}
\usepackage{latexsym}
\usepackage{amsmath}
\usepackage{amsthm}
\usepackage{amssymb}
\usepackage{mathrsfs}
\RequirePackage[colorlinks=true,citecolor=blue,urlcolor=blue]{hyperref}

\usepackage{caption}
\usepackage{subcaption}

\usepackage{eucal}

\usepackage{xcolor,soul}

\usepackage{mathtools}

 \usepackage[applemac]{inputenc}

\newtheorem{thm}{Theorem}[section]
\newtheorem{theo}[thm]{Theorem}

\newtheorem{lem}[thm]{Lemma}
\newtheorem{cor}[thm]{Corollary}
\newtheorem{defi}[thm]{Definition}
\newtheorem{hyp}[thm]{Assumption}

\theoremstyle{remark}
\newtheorem{ex}[thm]{Example}

\newtheorem{rem}[thm]{Remark}

\newcommand{\inP} {\underset{n \rightarrow \infty}{\overset{\mathbb{P}}{\xrightarrow{\hspace*{0.75cm}}}} }

\newcommand{\vip}{\vskip.2cm}
\newcommand{\COMMENTAIRE}[1]{}

\newcommand{\field}[1]{\mathbb{#1}}
\newcommand{\EE}{\field{E}}

\newcommand{\GG}{\field{G}}

\newcommand{\NN}{\field{N}}

\newcommand{\TT}{\field{T}}

\newcommand{\Hh}{{\mathcal H}}

\newcommand{\Pp}{{\mathcal P}}

\newcommand{\Qq}{{\mathcal Q}}

\def \ep {\varepsilon}

\newcommand{\rd}{{\rm d}}

\newcommand{\bff}{{\bf f}}
\newcommand{\bF}{{\mathfrak f}}

\newcommand{\cb}{{\mathcal B}}
\newcommand{\cc}{{\mathcal C}}

\newcommand{\ce}{{\mathcal E}}

\newcommand{\cf}{{\mathcal F}}

\newcommand{\ch}{{\mathcal H}}

\newcommand{\cn}{{\mathcal N}}

\newcommand{\cp}{{\mathcal P}}
\newcommand{\cq}{{\mathcal Q}}
\newcommand{\crr}{{\mathcal R}}

\newcommand{\cs}{{\mathscr S}}

\newcommand{\C}{{\mathbb C}}

\newcommand{\E}{{\mathbb E}}

\newcommand{\G}{\mathbb{G}}
\newcommand{\N}{{\mathbb N}}
\renewcommand{\P}{{\mathbb P}}
\newcommand{\p}{{\hat p}}

\newcommand{\R}{{\mathbb R}}

\newcommand{\T}{\mathbb{T}}

\newcommand{\ind}{{\bf 1}}

\newcommand{\sot}{\otimes_{\rm sym}}
\newcommand{\ssub}{\Sigma^{\rm sub}}

\newcommand{\scrit}{\Sigma^{\rm crit}}

\newcommand{\norm}[1]{\mathop{\parallel\! #1 \! \parallel}\nolimits}

\newcommand{\val}[1]{\mathop{\left| #1 \right|}\nolimits}
\newcommand{\inv}[1]{\mathop{\frac{1}{ #1}}\nolimits}

\newcommand{\reff}[1]{(\ref{#1})}

\begin{document}

\title[CLT for bifurcating {M}arkov chains under point-wise ergodic conditions]{Central limit theorem for bifurcating {M}arkov chains under point-wise ergodic conditions}

\author{S. Val\`ere Bitseki Penda and Jean-Fran\c cois Delmas}

\address{S. Val\`ere Bitseki Penda, IMB, CNRS-UMR 5584, Universit\'e Bourgogne Franche-Comt\'e, 9 avenue Alain Savary, 21078 Dijon Cedex, France.}

\email{simeon-valere.bitseki-penda@u-bourgogne.fr}

\address{Jean-Fran\c cois Delmas, CERMICS, Ecole des Ponts,  France.}

\email{jean-francois.delmas@enpc.fr}

\begin{abstract}
Bifurcating Markov  chains (BMC) are  Markov chains indexed by  a full binary tree representing  the evolution of a trait  along a population where each  individual has  two children.  We provide  a central limit theorem for  general additive functionals of BMC,  and prove the existence of three regimes. This  corresponds to a competition between the  reproducing  rate  (each  individual has  two  children)  and  the ergodicity rate  for the evolution of  the trait. This is  in contrast with  the  work  of  Guyon   (2007),  where  the  considered  additive functionals are  sums of  martingale increments,  and only  one regime appears.  Our result  can be seen as a discrete  time version, but with general trait evolution, of results in the time continuous setting of branching particle system from Adamczak and Mi\l{}o\'{s} (2015), where the evolution  of the trait is given by an Ornstein-Uhlenbeck process.
\end{abstract}

\maketitle






\section{Introduction}
Bifurcating Markov chains are a class of stochastic processes indexed by
regular binary tree and which satisfy the branching Markov property (see
below for a precise definition).  This model represents the evolution of
a trait along  a population where each individual has  two children.  To
the  best of  our knowledge,  the  term bifurcating  Markov chain  (BMC)
appears for the  first time in the work of  Basawa and Zhou \cite{BZ04}.
But,  it was  Guyon who,  in \cite{Guyon},  highlighted and  developed a
theory  of asymmetric  bifurcating  Markov chains.  Since  the works  of
Guyon, BMC theory  has been enriched from  probabilistic and statistical
point  of view  and several  extensions and  models using  BMC have  been
studied; we  can cite  the works  (see also  the references  therein) of
Bercu, de Saporta \&  G\'egout-Petit \cite{bdsgp:aa}, Delmas \& Marsalle
\cite{DM10},  Bitseki,   Djellout  \&  Guillin   \cite{BDG14},  Bitseki,
Hoffmann \& Olivier \cite{bho:aebmc},  Doumic, Hoffmann, Krell \& Robert
\cite{dhkr:segf}, Bitseki \& Olivier \cite{po:afe,bo:md} and Hoffmann \&
Marguet \cite{MR4025704}.  \medskip

The recent study of BMC models was motivated by the understanding of the
cell division mechanism (where the trait of an individual is given by its
growth rate).  The first model  of BMC, named  ``symmetric'' bifurcating
auto-regressive process  (BAR) 
were introduced  by Cowan  \& Staudte
\cite{CS86} in  order to analyze cell  lineage data. Since the  works of
Cowan  and Staudte,  many  extensions  of their  model  were studied  in
Markovian  and non-Markovian  setting (see  for e.g.   \cite{po:afe} and
references therein).  In particular,  in \cite{Guyon}, Guyon has studied
``asymmetric'' BAR 
in order to prove statistical  evidence  of aging  in  Escherichia
Coli, giving a new approach to the problem studied in
\cite{smpt:adormsd}. 
 Let us also note that BMC have been used recently
in several statistical works to study the estimator of the cell division
rate \cite{dhkr:segf,bho:aebmc,MR4025704}. Moreover, another studies,
such as  \cite{MR3906858}, can be generalized using the BMC
theory (we refer to the conclusion therein).
\medskip

In  this paper,  our objective  is to establish  a central  limit
theorem for  additive functionals of  BMC. 
With respect  to this objective,
notice that asymptotic results for BMC have been studied in \cite{Guyon}
(law of  large numbers  and central limit  theorem) and  in \cite{BDG14}
(moderate deviations  principle and  strong law  of large  numbers). See
\cite{DM10} for the  law of large numbers and central  limit theorem for
BMC on Galton-Watson  tree.  Notice also that  recently, limit theorems,
in  particular law  of large  numbers,  has been  studied for  branching
Markov process,  see \cite{MR4011569} and \cite{c:ltbm},  and that large
values  of  parameters in  stable  BAR  process  allows to  exhibit  two
regimes, see  \cite{BB2020}.  However, the central  limit theorems which
appear  in  \cite{Guyon,bdsgp:aa,DM10}  have   been  done  for  additive
functionals using increments of  martingale, which implies in particular
that the functions considered depend on the traits of the mother and its
two daughters.  The study of the case where the functions depend only on
the trait of  a single individual has  not yet been treated  for BMC (in
this case it  is not useful to  solve the Poisson equation  and to write
additive functional as  sums of martingale increments as  the error term
on  the  last  generation  is  not  negligible  in  general).  For  such
functions, the  central limit  theorems have  been studied  recently for
branching      Markov     processes      and     for      superprocesses
\cite{am:cltou,m:sclt,MR3592799,MR3803921}.
Our results  can be seen as  a discrete version of  those given in
the previous works, but with general ergodic hypothesis on the evolution
of the  trait.  Unlike the results  given in \cite{Guyon,bdsgp:aa,DM10},
we observe  three regimes  (sub-critical, critical  and super-critical),
which correspond to  a competition between the reproducing  rate (here a
mother has two  daughters) and the ergodicity rate for  the evolution of
the trait  along a lineage  taken uniformly at random.   This phenomenon
already appears  in the  works of Athreya  \cite{athreya1969limit}.  For
BMC  models,  we stress  that  the  three  regimes already  appears  for
moderate   deviations  and   deviation   inequalities  in   \cite{BDG14,
  bitsekidjellout2014,bitseki2015}.
\medskip

We follow the approach of    \cite{Guyon, DM10}  and consider
ergodic  theorem with respect to the 
point-wise  convergence.
However, unlike the latter papers, we provide  a different  normalization for  the fluctuations according to the regime being critical, sub-critical and super-critical,
see  respectively Corollaries  \ref{cor:subcritical}, \ref{cor:critical}
and  \ref{cor:super-crit-d1-2}.   We  shall explicit  in  a  forthcoming
paper, that  those results allow to  recover the one regime  result from
\cite{Guyon}  for additive  functionals  given by  a  sum of  martingale
increments.  \medskip

The  paper is  organized  as follows.   We introduce  the  BMC model  in
Section \ref{sec:BMC-sub}  and consider the  sets of assumptions  in the
spirit  of  \cite{Guyon}  in Section~\ref{sec:pw-approach}.   The  main
results    are     presented    in     Section~\ref{sec:results}:    see
Section~\ref{sec:main-res-sub}  for results  in  the sub-critical  case,
with    technical    proofs    in    Section~\ref{sec:proof-sub};    see
Section~\ref{sec:main-res-crit} for  results in the critical  case, with
technical    proofs    in    Section~\ref{sec:proof-crit};    and    see
Section~\ref{sec:main-res-super} for results in the super-critical case,
with  technical  proofs  in  Section~\ref{sec:proof-super}.   The  proof
relies essentially  on explicit second moments  computations and precise
upper  bounds  of  fourth  moments   for  BMC,  which  are  recalled  in
Section~\ref{sec:moment-BMC}.

\section{Models and assumptions}
\label{sec:BMC}
\subsection{Bifurcating Markov chain: the model}
\label{sec:BMC-sub}

We denote by   $\N$     the   set   of  non-negative   integers   and
$\N^*=   \N   \setminus   \{0\}$.
If $(E,  \ce)$ is a  measurable space, then $\cb(E)$  (resp. $\cb_b(E)$,
resp.    $\cb_+(E)$)  denotes   the  set   of  (resp.    bounded,  resp.
non-negative)  $\R$-valued measurable  functions  defined  on $E$.   For
$f\in \cb(E)$, we set $\norm{f}_\infty =\sup\{|f(x)|, \, x\in E\}$.  For
a finite measure $\lambda$ on $(E,\ce)$ and $f\in \cb(E)$ we shall write
$\langle \lambda,f  \rangle$ for  $\int f(x) \,  \rd\lambda(x)$ whenever
this  integral is  well defined.  For $n\in \N^*$, the product space $E^n$ is
endowed with the product $\sigma$-field  $\ce^{\otimes n}$. 
If $(E,d)$  is a  metric space,  then
$\ce $ will denote its Borel $\sigma$-field and the set $\cc_b(E)$ (resp.
$\cc_+(E)$) denotes the set of bounded (resp.  non-negative) $\R$-valued
continuous functions defined on $E$. 
\medskip

Let  $(S, \cs)$   be  a  measurable  space. 
Let $Q$ be a   
probability kernel   on $S \times \cs$, that is:
$Q(\cdot  , A)$  is measurable  for all  $A\in \cs$,  and $Q(x,
\cdot)$ is  a probability measure on $(S,\cs)$ for all $x \in
S$. For any $f\in \cb_b(S)$,   we set for $x\in S$:
\begin{equation}
   \label{eq:Qf}
(Qf)(x)=\int_{S} f(y)\; Q(x,\rd y).
\end{equation}
We define $(Qf)$, or simply $Qf$, for $f\in \cb(S)$ as soon as the
integral \reff{eq:Qf} is well defined, and we have $\cq f\in \cb(S)$. For $n\in \N$, we denote by $Q^n$  the
$n$-th iterate of $Q$ defined by $Q^0=I_d$, the identity map on
$\cb(S)$, and $Q^{n+1}f=Q^n(Qf)$ for $f\in \cb_b(S)$.  

Let $P$ be a   
probability kernel   on $S \times \cs^{\otimes 2}$, that is:
$P(\cdot  , A)$  is measurable  for all  $A\in \cs^{\otimes 2}$,  and $P(x,
\cdot)$ is  a probability measure on $(S^2,\cs^{\otimes 2})$ for all $x \in
S$. For any $g\in \cb_b(S^3)$ and $h\in \cb_b(S^2)$,   we set for $x\in S$:
\begin{equation}
   \label{eq:Pg}
(Pg)(x)=\int_{S^2} g(x,y,z)\; P(x,\rd y,\rd z)
\quad\text{and}\quad
(Ph)(x)=\int_{S^2} h(y,z)\; P(x,\rd y,\rd z).
\end{equation}
We define $(Pg)$ (resp. $(Ph)$), or simply $Pg$ for $g\in \cb(S^3)$
(resp. $Ph$ for $h\in \cb(S^2)$), as soon as the corresponding 
integral \reff{eq:Pg} is well defined, and we have  that $Pg$ and
$Ph$ belong to $\cb(S)$.
\medskip 

We  now introduce  some notations  related to  the regular  binary tree.
 We   set   $\T_0=\G_0=\{\emptyset\}$,
$\G_k=\{0,1\}^k$  and $\T_k  =  \bigcup _{0  \leq r  \leq  k} \G_r$  for
$k\in  \N^*$, and  $\T  =  \bigcup _{r\in  \N}  \G_r$.   The set  $\G_k$
corresponds to the  $k$-th generation, $\T_k$ to the tree  up to the $k$-th
generation, and $\T$ the complete binary  tree. For $i\in \T$, we denote
by $|i|$ the generation of $i$ ($|i|=k$  if and only if $i\in \G_k$) and
$iA=\{ij; j\in A\}$  for $A\subset \T$, where $ij$  is the concatenation
of   the  two   sequences  $i,j\in   \T$,  with   the  convention   that
$\emptyset i=i\emptyset=i$.

We recall the definition of bifurcating Markov chain  from
\cite{Guyon}. 
\begin{defi}
  We say  a stochastic process indexed  by $\T$, $X=(X_i,  i\in \T)$, is
  a bifurcating Markov chain (BMC) on a measurable space $(S, \cs)$ with
  initial probability distribution  $\nu$ on $(S, \cs)$ and probability
  kernel $\cp$ on $S\times \cs^{\otimes 2}$ if:
\begin{itemize}
\item[-] (Initial  distribution.) The  random variable  $X_\emptyset$ is
  distributed as $\nu$.
   \item[-] (Branching Markov property.) For  a sequence   $(g_i, i\in
     \T)$ of functions belonging to $\cb_b(S^3)$, we have for all $k\geq 0$,
\[
\E\Big[\prod_{i\in \G_k} g_i(X_i,X_{i0},X_{i1}) |\sigma(X_j; j\in \T_k)\Big] 
=\prod_{i\in \G_k} \cp g_i(X_{i}).
\]
\end{itemize}
\end{defi}

Let $X=(X_i,  i\in \T)$ be a BMC  on a measurable space $(S, \cs)$ with
  initial probability distribution  $\nu$ and probability
  kernel $\cp$. 
We    define     three
probability kernels $P_0, P_1$ and $\cq$ on $S\times \cs$ by:
\[
P_0(x,A)=\cp(x, A\times S), \quad
P_1(x,A)=\cp(x, S\times A) \quad\text{for
$(x,A)\in S\times  \cs$, and}\quad
\cq=\inv{2}(P_0+P_1).
\] 
Notice  that  $P_0$ (resp.   $P_1$)  is  the  restriction of  the  first
(resp. second) marginal of $\cp$ to $S$.  Following \cite{Guyon}, we
introduce an  auxiliary Markov  chain $Y=(Y_n, n\in  \N) $  on $(S,\cs)$
with  $Y_0$ distributed  as $X_\emptyset$  and transition  kernel $\cq$.
The  distribution of  $Y_n$ corresponds  to the  distribution of  $X_I$,
where $I$  is chosen independently from  $X$ and uniformly at  random in
generation  $\G_n$.    We  shall   write  $\E_x$   when  $X_\emptyset=x$
(\textit{i.e.}  the initial  distribution  $\nu$ is  the  Dirac mass  at
$x\in S$).  \medskip

We  end  this  section  with  a  useful  notation.  By  convention,  for
$f,g\in  \cb(S)$, we  define the  function $f\otimes  g\in \cb(S^2)$  by
$(f\otimes g)(x,y)=f(x)g(y)$ for $x,y\in S$ and introduce the notations:
\[
f\sot g= \inv{2}(f\otimes g + g\otimes f) \quad\text{and}\quad f\otimes ^2= f\otimes f.
\]
Notice that 
$\cp(g\sot \ind)=\cq(g)$ for  $g\in \cb_+(S)$.

\subsection{Assumptions}
\label{sec:pw-approach}

For  a  set  $F\subset  \cb(S)$   of  $\R$-valued  functions,  we  write
$F^2=\{f^2; f\in F\} $, $F\otimes  F=\{f_0\otimes f_1; f_0, f_1\in F\}$,
and  $P(E)=\{Pf;  f\in E\}$  whenever  a  kernel $P$  act  on  a set  of
functions $E$.  Following \cite{Guyon}, we state a structural assumption
on the set of functions we shall consider.

\begin{hyp}
   \label{hyp:F}
Let $F\subset \cb(S)$ be a set of $\R$-valued functions such that:
\begin{itemize}
   \item[$(i)$] $F$ is a  vector subspace which contains the constants;
   \item[$(ii)$] $F^2 \subset F$;
   \item[$(iii)$] $F\subset L^1(\nu)$; 
 \item[$(iv)$]  $F\otimes F \subset L^1(\cp(x, \cdot))$ for all $x\in S$,
    and $\cp(F\otimes F)\subset F$.
\end{itemize}
\end{hyp}

The   condition   $(iv)$   implies   that   $P_0(F)\subset   F$,
$P_1(F)\subset F$ as  well as $\cq(F)\subset F$.  Notice that if  $f\in
F$, then even
if $|f|$ does  not belong to $F$, using  conditions  $(i)$  and $(ii)$,
we get, with $g=(1+f^2)/2$, that   $|f|\leq  g$ and $g\in F$. Typically,
when $(S, d)$ is a metric space, the set $F$  can be the set  $\cc_b(S)$
of bounded real-valued functions, or the set of smooth real-valued functions such
that all derivatives have  at
most polynomials growth. 
\medskip

Following \cite{Guyon}, we also  consider the following  ergodic properties for $\cq$. 
\begin{hyp}
   \label{hyp:F1}
There exists a probability measure $\mu$ on $(S, \cs)$ 
  such that $F\subset L^1(\mu)$ and for all $f\in F$, we have the
  point-wise convergence  $\lim_{n\rightarrow \infty } \cq^{n}f =
  \langle \mu, f \rangle$ and 
there exists $g\in F$ with:
\begin{equation}
   \label{eq:erg-bd}
 |\cq^n(f)|\leq  g\quad\text{for all $n\in
  \N$.}
 \end{equation} 
\end{hyp}

We consider also the following geometrical ergodicity. 
\begin{hyp}
   \label{hyp:F2}
There exists a probability measure $\mu$ on $(S, \cs)$  such that
$F\subset L^1(\mu)$,  and $\alpha\in (0, 1)$ such that for all $f\in F$
  there exists $g\in F$ such that: 
 \begin{equation}
   \label{eq:geom-erg}
|\cq^{n}f - \langle \mu, f \rangle| \leq \alpha^{n} g \quad
\text{for all  $n\in \N$.}
\end{equation}
\end{hyp}

A  sequence $\bF=(f_\ell,  \ell\in  \N)$ of  elements  of $F$  satisfies
uniformly \reff{eq:erg-bd}  and \reff{eq:geom-erg} if there  is $g\in F$
such that:
\begin{equation}\label{eq:unif-f}
|\cq^n(f_\ell)|\leq  g \quad\text{and}\quad |\cq^{n}f_\ell  - \langle \mu, f_\ell \rangle|  \leq \alpha^{n} g \quad \text{for all $n,\ell\in \N$.}
\end{equation}
This    implies    in    particular   that    $|f_\ell|\leq    g$    and
$|\langle \mu, f_\ell \rangle|\leq \langle \mu, g \rangle$.  Notice that
\reff{eq:unif-f} trivially holds if $\bF$ takes finitely distinct values
(\textit{i.e.} the subset $\{f_\ell; \ell\in  \N\}$ of $F$ is finite)
each satisfying \reff{eq:erg-bd}  and \reff{eq:geom-erg}. 

\begin{ex}
   \label{ex:metric-space}
Let  $(S, d)$ be a metric space, $\cs$ its Borel $\sigma$-field, and
  $Y$  a   Markov chain  uniformly
   geometrically ergodic \emph{i.e.}  there exists $\alpha \in (0,1)$
   and  a  finite  constant  $C$  such  that  for all  $x \in S$:
\begin{equation}\label{eq:unif-erg}
\|\Qq^{n}(x,\cdot) - \mu\|_{TV} \leq C \alpha^{n},
\end{equation}
where,  for  a signed  finite  measure  $\pi$  on $(S,\cs)$,  its  total
variation            norm            is            defined            by
$\|\pi\|_{TV}   =     \sup_{f\in \cb(S), \, \norm{f}_\infty  \leq   1}   |\langle
\pi,f\rangle|$.
Then,  taking for $F$  the  set of  $\R$-valued continuous  bounded
function   $\cc_b(S)$,  we   get   that  properties (i-iii) from
Assumption  \ref{hyp:F}   and Assumption 
\ref{hyp:F2} hold. In particular, Equation  \eqref{eq:unif-erg} implies
that 
\reff{eq:geom-erg} holds with $g=C \norm{f}_{\infty}$. 
\end{ex}

\medskip

We consider the stronger ergodic property based on a second spectral
gap. 
\begin{hyp}
   \label{hyp:F3}
There exists a probability measure $\mu$ on $(S, \cs)$  such that
$F\subset L^1(\mu)$,  and $\alpha\in (0, 1)$, a finite  non-empty set
  $J$ of indices, distinct complex eigenvalues $\{\alpha_j, \, j\in J\}$
  of  the  operator  $\cq$ with  $|\alpha_j|=\alpha$,  non-zero  complex
  projectors $\{\crr_j, \, j\in J\}$  defined on $\C F$, the $\C$-vector
  space    spanned by        $F$,       such        that
  $\crr_j\circ \crr_{j'}=\crr_{j'}\circ  \crr_{j}=0$ for all  $j\neq j'$
  (so that  $\sum_{j\in J} \crr_j$ is  also a projector defined  on $\C F$)
  and a positive  sequence $(\beta_n, n\in \N)$ converging  to $0$, such
  that  for  all  $f\in  F$  there  exists  $g\in  F$  and,  with
  $\theta_j=\alpha_j/\alpha$:
\begin{equation}
   \label{eq:hyp-crit}
\Big|\cq^{n}( f) - \langle \mu, f \rangle - \alpha^n \sum_{j\in J}
\theta_j^n\,  \crr_j (f) \Big| \leq 
  \beta_n \alpha^{n} g \quad
\text{for all  $n\in \N$.}
\end{equation}
\end{hyp}

Without loss of generality, we shall assume that the sequence $(\beta_n,
n\in \N)$ in Assumption  \ref{hyp:F3} is non-increasing and bounded from above by 1.

\begin{rem}
   \label{rem:Guyon} 
   In  \cite{Guyon},  only  the   structural Assumption \ref{hyp:F}  and
   the ergodic  Assumption  
   \ref{hyp:F1}  were assumed.  If $F$ contains a set $A$ of bounded
   functions which is separating (that is two probability measures which
   coincides on $A$ are equal), then Assumption  \ref{hyp:F}  and
    \ref{hyp:F1}
   imply in particular that $\mu$ is
   the only invariant measure of $\cq$.  Notice  that  the geometric  ergodicity
   Assumption  \ref{hyp:F2} implies  Assumption  \ref{hyp:F1}, and  that
   Assumption  \ref{hyp:F3} implies  Assumption  \ref{hyp:F2} (with  the
   same $\alpha$ but possibly different function $g$).  \medskip
\end{rem}

\begin{ex}
  \label{ex:sBAR}
We  consider the real-valued Gaussian symmetric bifurcating autoregressive
process (BAR) $X=(X_{u},u\in\TT)$ where for all $ u \in \T\backslash \{\emptyset\}$: 
\[
  X_{u} = a X_{v} + \ep_{v},
\]
where $v$ is the parent of $u$, that is $u=v0$ or $u=v1$, $a\in (-1,
1)$,  and $(\ep_v, \, v\in \TT)$ are independent  Gaussian
random variables $\cn(0, \sigma^2)$ with  $\sigma>0$. We obtain:
\[
  \cp(x, dy, dz)=\cq(x, dy) \cq(x, dz)
  \quad\text{with}\quad \cq f(x)=\E[f(ax +\sigma G)],
\]
where $G$ is a standard $\cn(0, 1)$ Gaussian random variable.
More generally we have $\cq^n
f(x)=\E\left[f\left(a^n x + \sqrt{1- a^{2n}} \sigma_a G\right)\right]$,
where 
$\sigma_a=\sigma (1- a^2)^{-1/2}$. 
The kernel 
$\Qq$ admits a unique invariant probability measure $\mu$, which is Gaussian
$\cn(0, \sigma_a^2)$. The    operator $\Qq$ (on
$L^2(\mu)$) is a symmetric integral Hilbert-Schmidt operator whose eigenvalues are given
by $\sigma_{p}(\Qq) = (a^{n}, n \in \NN)$, their algebraic
multiplicity is one and the corresponding eigen-functions 
$(\bar{g}_{n}(x), n \in \NN)$ are defined for $n\in \N$ by $
\bar{g}_{n}(x) = g_{n}\left(\sigma_a^{-1} \, x\right)$,
where $g_{n}$ is the Hermite polynomial of degree $n$. In particular, we
have $\bar g_0=1$ and
$\bar g_1(x)=\sigma_a^{-1} x$. Let  $\crr$ be the  orthogonal projection on the vector
space generated by $\bar{g}_{1}$, that is $
\crr f= \langle \mu, f\bar g_1 \rangle\, \bar g_1$ or equivalently, for
$x\in \R$:
\begin{equation}
   \label{eq:R-symBAR}
 \crr f(x)=\sigma_a^{-1}\, x \, \E\left[G f(\sigma_aG)\right].
\end{equation}

Consider $F$ the set of functions $f\in \cc^2(\R)$ such that $f, f'$ and
$f''$ have at most polynomial growth. 
And assume that the
probability distribution $\nu$ has all its moments, which is equivalent
to say that $F\subset L^1(\nu)$. Then the set $F$ satisfies Assumption
\ref{hyp:F}. We also have that $F\subset L^1(\mu)$. Then, it is not
difficult to check directly that Assumption 
\ref{hyp:F3} also holds with $J=\{j_0\}$, $\alpha_{j_0} = \alpha = a$,
$\beta_n =  a^{n}$ and $\crr_{j_0}=\crr$ (and also  Assumptions \ref{hyp:F1} and \ref{hyp:F2}
hold).
\end{ex}

\subsection{Notations for average of different functions over different generations}
\label{sec:more-not}

Let  $X=(X_u, u\in  \T)$ be  a BMC  on $(S, \cs)$ with  initial
 probability distribution $\nu$, and probability kernel $\cp$.  Recall
   $\cq$ is the induced Markov kernel. We  assume that $\mu$ is an invariant probability measure of $\cq$. 
   
For a finite set $A\subset \T$ and a function $f\in \cb(S)$, we set:
\[
M_A(f)=\sum_{i\in A} f(X_i).
\]
We shall be interested in the cases $A=\G_n$ (the  $n$-th generation)
and $A=\T_n$ (the tree up to the $n$-th generation). We recall  from \cite[Theorem~11 and Corollary~15]{Guyon}  that under Assumptions
\ref{hyp:F} and \ref{hyp:F1} (resp. and also Assumption
\ref{hyp:F2}), we have for $f\in F$ the following convergence in $L^2(\mu)$
(resp. a.s.):
 \begin{equation}
    \label{eq:lfgn-G}
 \lim_{n\rightarrow\infty } |\G_n|^{-1} M_{\G_n}(f)=\langle \mu, f
 \rangle
 \quad\text{and}\quad
 \lim_{n\rightarrow\infty } |\T_n|^{-1} M_{\T_n}(f)=\langle \mu, f
 \rangle.
 \end{equation}
 
We shall now consider the corresponding fluctuations.  We will use
frequently the following notation: 
\[
\boxed{\tilde f= f - \langle \mu,
    f \rangle\quad \text{for $f\in L^1(\mu)$.}}
\]

In order to study the asymptotics of $ M_{\GG_{n-\ell }}(\tilde
f)$,  we shall consider the contribution of the descendants of the
individual $i\in \T_{n-\ell}$ for  $n\geq \ell\geq 0$:
\begin{equation}
   \label{eq:def-Nnil}
N^\ell_{n,i}(f)=|\G_n|^{-1/2} M_{i\G_{n-|i|-\ell}}(\tilde f), 
\end{equation}
where  $i\G_{n-|i|-\ell}=\{ij, \, j\in
\G_{n-|i|-\ell}\}\subset \G_{n-\ell}$. For  all $k\in \N$ such that $n\geq k+\ell$, we have:
\[
M_{\G_{n-\ell}}(\tilde f)=\sqrt{|\G_n|}\,\, \sum_{i\in \G_k}
N^\ell_{n,i}(f)=
\sqrt{|\G_n|}\, \,N_{n, \emptyset}^\ell(f).
\]
Let $\bF=(f_\ell, \ell\in \N)$ be a sequence of elements of
$L^1(\mu)$. We set
for $n\in \N$ and $i\in \T_n$:
\begin{equation}
   \label{eq:def-NiF}
N_{n,i}(\bF)=\sum_{\ell=0}^{n-|i|} N_{n,i}^\ell(f_\ell) 
=|\G_n|^{-1/2 }\sum_{\ell=0}^{n-|i|}
M_{i\G_{n-|i|-\ell}}(\tilde f_\ell).
\end{equation}
In $N_{n,  i}$,   we consider  the  
contribution of the  descendants of $i$ up to generation  $n$. We deduce
that
$   \sum_{i\in   \G_k}  N_{n,i}(\bF)=|\G_n|^{-1/2   }\sum_{\ell=0}^{n-k}
M_{\G_{n-\ell}}(\tilde f_\ell)$ which gives for $k=0$:
\begin{equation}
   \label{eq:def-NOf}
\boxed{N_{n, \emptyset}(\bF)= |\G_n|^{-1/2 }\sum_{\ell=0}^{n}
  M_{\G_{n-\ell}}(\tilde f_{\ell}).}
\end{equation}
In  $N_{n,  \emptyset}$,  we  consider   the  contribution  of  all  the
individual from generation $0$ up to  generation $n$. We shall prove the
convergence in law of $N_{n, \emptyset}(\bF)$ in the following sections.

\begin{rem}
   \label{rem:simpleN0n}
 We shall consider in particular the following two simple cases. 
Let      $f\in     L^1(\mu)$      and     consider      the     sequence
$\bF=(f_\ell,  \ell\in \N)$.   If  $f_0=f$ and $f_\ell=0$ for  $\ell
\in \N^*$,  then we get:
\[
N_{n, \emptyset}(\bF)= |\G_n|^{-1/2} M_{\G_n}(\tilde f).
\]
  If $f_\ell=f$ for $\ell \in \N$, then we shall write  $\bff=(f,f,
  \ldots)$, and we get, as  $|\T_n|=2^{n+1} - 1 $ and $|\G_n|= 2^n$:
\[
N_{n, \emptyset}(\bff)= |\G_n|^{-1/2} M_{\T_n}(\tilde f)
=  \sqrt{2 - 2^{-n}}\,\, |\T_n|^{-1/2} M_{\T_n}(\tilde f).
\]
Thus, we will easily deduce the
fluctuations of $M_{\T_n}(f)$ and $M_{\G_n}(f)$ from the asymptotics of $N_{n,
  \emptyset}(\bF)$. 
\end{rem}

\medskip

To study the asymptotics of $N_{n, \emptyset}(\bF)$, it is convenient to
write for $n\geq k\geq 1$:
\begin{equation}
   \label{eq:nof-D}
N_{n, \emptyset}(\bF)= |\G_n|^{-1/2} \sum_{r=0}^{k-1} M_{\G_r}(\tilde
f_{n-r}) + \sum_{i\in \G_k} N_{n,i}(\bF).
\end{equation}
If $\bff=(f,f, \ldots)$ is the infinite sequence of the same function
$f$, this becomes:
\begin{equation}
   \label{eq:Nnff}
N_{n, \emptyset}(\bff)= |\G_n|^{-1/2} M_{\T_n}(\tilde f)= |\G_n|^{-1/2}
M_{\T_{k-1}}(\tilde f)+ \ \sum_{i\in
  \G_k} N_{n,i}(\bff).
\end{equation}
\medskip

In the proofs, we will  denote by  $C$ any
unimportant finite  constant which may  vary from  line to line  (in
particular $C$ does not  depend on $n\in \N$ nor on the considered sequence of
functions $\bF=(f_\ell, \ell \in \N)$).

\section{Main results} \label{sec:results}

\subsection{The sub-critical case: $2\alpha^2<1$}\label{sec:main-res-sub}

We shall consider, when well defined, for a sequence $\bF=(f_\ell,
\ell\in  \N)$ of measurable real-valued
functions defined on $S$, the quantities:
\begin{equation}\label{eq:ssub}
\ssub(\bF)=\ssub_1(\bF)+ 2 \ssub_2(\bF), 
\end{equation}
where:
\begin{align}
   \label{eq:S1}
   \ssub_1(\bF)
&=\sum_{\ell\geq 0} 
2^{-\ell} \, \langle \mu,   \tilde f_\ell^ 2\rangle
+ \sum_{\ell\geq 0, \, k\geq 0} 
2^{k-\ell} \, \langle \mu, \cp\left((\cq^k \tilde
    f_\ell) \otimes^2\right)\rangle,\\
   \label{eq:S2}
   \ssub_2(\bF)
&=\sum_{0\leq \ell< k}
2^{-\ell} \langle \mu,   \tilde f_k \cq^{k-\ell} \tilde f_\ell\rangle
+ \sum_{\substack{0\leq \ell< k\\ r\geq 0}} 2^{r-\ell} \langle \mu,
\cp\left( \cq^r \tilde
    f_k \sot
\cq^{k-\ell+r} \tilde f_\ell  \right)\rangle.
\end{align}
\medskip

We have the following result whose proof is given in Section
\ref{sec:proof-sub}. 
\begin{theo}\label{theo:subcritical}
Let  $X$  be  a  BMC   with  kernel  $\cp$  and  initial  distribution
$\nu$ such that Assumptions  \ref{hyp:F} and  \ref{hyp:F2}  are  in force  with $\alpha\in  (0, 1/\sqrt{2})$.  We  have the  following convergence  in distribution for  all sequence $\bF=(f_\ell, \ell\in  \N)$ of elements of $F$ satisfying Assumptions  \ref{hyp:F2} uniformly, that is \reff{eq:unif-f} for some $g\in F$:
\begin{equation*}
N_{n, \emptyset}(\bF) \; \xrightarrow[n\rightarrow \infty ]{\text{(d)}} \; G,
\end{equation*}
where $G$ is a centered Gaussian  random variable with  variance
$\ssub(\bF)$ given  by \reff{eq:ssub}, which is well defined and
finite. 
\end{theo}

The  convergence in  distribution of  $N_{n, \emptyset}(\bF)$  allows to
recover the  convergence in distribution  of the average  over different
successive                                                   generations
$|\GG_{n}|^{-1/2}  (M_{\GG_{n   }}(\tilde  f_{0}),\ldots,   M_{\GG_{n  -
    k}}(\tilde f_{k}))$.  Notice the limit  is a Gaussian  random vector
$(G_1,   \ldots,   G_k)$.   \emph{A   priori}   the   random   variables
$G_1, \ldots,  G_k$ are  not independent because of the interaction
coming from 
\eqref{eq:S2}.  In contrast,  it was  proved in  \cite{DM10}, that  the
average over  different successive generations of  martingale increments
converges to Gaussian independent random variables.

\begin{rem}\label{rem:sub-cr-F=ff}
For $f\in \cb(S)$, when it is well defined, we set:
\begin{equation}
\label{eq:ST} 
\ssub_{\GG}(f) = \langle \mu,   \tilde f^ 2\rangle + \sum_{k\geq 0}
2^{k} \, \langle \mu, \cp\left(\cq^k \tilde f \otimes ^2\right)\rangle
\quad\text{and}\quad
\ssub_{\TT}(f) = \ssub_{\G}(f) +2  \ssub_{\TT,2}(f),
\end{equation}
where
\[
  \ssub_{\TT,2}(f) 
=    \sum_{k\geq 1} \langle \mu, \tilde  f \cq^{k}
\tilde  f\rangle
+  \sum_{\substack{k\geq 1 \\   r\geq  0}}
    2^{r}   \langle   \mu,   \cp\left(  \cq^r   \tilde   f   \sot
      \cq^{r+k} \tilde f \right)\rangle.
\]

If we take $\bF = (f,0,0, \ldots)$, we have $\ssub(\bF)=\ssub_{\G}(f)$. If we take $ \bff=(f,f, \ldots)$, the infinite sequence of the same function $f$, we have $\ssub(\bff)=2\ssub_{\TT}(f) $.
\end{rem}

As a direct consequence of Remarks \ref{rem:sub-cr-F=ff} and
\ref{rem:simpleN0n}, and the more general Theorem
\ref{theo:subcritical}, we get
the following result. 

\begin{cor}\label{cor:subcritical}
  Let $X$ be a BMC with kernel $\cp$ and initial distribution $\nu$ such
  that     Assumptions     \ref{hyp:F}    and     \ref{hyp:F2} are in  force     with
  $\alpha\in (0,  1/\sqrt{2})$.  Let $f\in F$.   Then, we
  have the following convergence in distribution:
\begin{equation*}
|\GG_{n}|^{-1/2}M_{\GG_{n}}(\tilde{f}) \; \xrightarrow[n\rightarrow \infty ]{\text{(d)}} \; G_{1}
  \quad \text{and} \quad
  |\TT_{n}|^{-1/2}M_{\TT_{n}}(\tilde{f}) \; \xrightarrow[n\rightarrow
  \infty ]{\text{(d)}} \; G_{2},  
\end{equation*}
where $G_{1}$ and $G_2$ are  centered Gaussian  random variables with
respective  
variances   $\ssub_{\GG}(f)$ and $\ssub_{\TT}(f)$ 
given in \reff{eq:ST},   which are  well defined
and finite. 
\end{cor}
 
\begin{proof}[Proof of Corollary \ref{cor:subcritical}]
Take
the infinite sequence $\bF = (f,0,0,\cdots)$, where only the first
component is non-zero, to deduce from Theorem \ref{theo:subcritical} the
convergence in distribution of
$|\G_{n}|^{-1/2}M_{\GG_{n}}(\tilde{f})=N_{n, \emptyset} (\bF)  $.
Next, take the infinite sequence $\bF = \bff=(f,f, \ldots)$ of the same function
$f$ in  Theorem \ref{theo:subcritical} and use \reff{eq:Nnff} and as
well as 
$\lim_{n\rightarrow \infty} |\GG_{n}|/|\TT_{n}| = 1/2$, to get  the
convergence in distribution  for $|\TT_{n}|^{-1/2}M_{\TT_{n}}(\tilde{f})
=(|\G_n|/|\T_n|)^{1/2} N_{n, \emptyset} (\bff)$. 
\end{proof}

 

\subsection{The critical case $2\alpha^2=1$}
\label{sec:main-res-crit}

In the  critical case $\alpha=1/\sqrt{2}$,  we shall denote  by $\crr_j$
the projector on the eigen-space associated to the eigenvalue $\alpha_j$
with  $\alpha_j=\theta_j  \alpha$, $|\theta_j|=1$  and  for  $j$ in  the
finite set of indices  $J$. Since $\cq$ is a real  operator, we get that
if $\alpha_j$ is a non real  eigenvalue, so is $\overline \alpha_j$.  We
shall  denote   by  $\overline  \crr_j$  the   projector  associated  to
$\overline \alpha_j$.  Recall that the  sequence $(\beta_n, n\in \N)$ in
Assumption  \ref{hyp:F3} can  (and  will) be  chosen non-increasing  and
bounded from  above by  1. For all  measurable real-valued  function $f$
defined on $S$, we set, when this is well defined:
\begin{equation}\label{eq:fhatcrit-S}
\boxed{ \hat{f} = \tilde{f} - \sum_{j \in J} \crr_{j}(f)
\quad\text{with}\quad
\tilde f= f- \langle \mu, f \rangle.}
\end{equation}
We shall consider, when well defined, for a sequence $\bF=(f_\ell,
\ell\in  \N)$ of measurable real-valued
functions defined on $S$, the quantities:
\begin{equation}\label{eq:scrit}
\scrit (\bF)=\scrit _1(\bF)+ 2\scrit _2(\bF), 
\end{equation}
where:
\begin{align}
\label{eq:S1-crit}
\scrit_1(\bF) 
&= \sum_{k\geq 0} 2^{-k} \langle \mu, \cp f_{k, k}^*
\rangle
= \sum_{k\geq 0} 2^{-k} \sum_{j\in J} \langle \mu,
\cp(\crr_{j}(f_k) \sot \overline{\crr}_{j}(f_k)) \rangle,\\
\label{eq:S2-crit}
\scrit_2(\bF) 
&=   \sum_{0\leq \ell <k} 2^{-(k+\ell)/2} \langle \mu, \cp
f_{k, \ell}^*  \rangle,
\end{align}
with, for $k, \ell\in \N$:
\begin{equation}\label{eq:def-fkl*}
f_{k, \ell}^*= \sum_{j\in J}\,\, \theta_j^{\ell -k}\,  \crr_j
(f_k)\sot   \overline \crr_j (f_\ell).  
\end{equation}
Notice  that $f_{k, \ell}^*=f_{\ell, k}^*$ and that $f_{k, \ell}^*$ is
real-valued as  $\overline {\theta_j^{\ell -k}\,  \crr_j (f_k)\otimes
\overline \crr_j (f_\ell)}= \theta_{j'}^{\ell -k}\,  \crr_{j'} (f_k)\otimes
\overline \crr_{j'} (f_\ell)$ for $j'$ such that $\alpha_{j'}=\overline
\alpha _j$ and thus $\crr_{j'}=\overline \crr_j$. 
\medskip

We shall  consider sequences $\bF=(f_\ell,  \ell\in \N)$ of  elements of
$F$ which satisfies Assumption \ref{hyp:F3} uniformly, that is such that
there exists $g\in F$ with:
\begin{equation}\label{eq:unif-f-crit}
|\cq^n(f_\ell)|\leq  g, \quad |\cq^{n}(\tilde f_\ell)| \leq \alpha^{n} g \quad\text{and}\quad |\cq^{n}(\hat{f}_\ell)| \leq \beta_n \, \alpha^{n} g \quad\text{for all $n,\ell\in \N$.}
\end{equation}
We deduce  that there exists a  finite constant $c_J$ depending  only on $\{\alpha_j, \,  j\in J\} $ such  that for all $\ell\in  \N$, $n\in \N$, $j_0\in J$:
\begin{equation}\label{eq:cons-f-crit}
|f_\ell|\leq    g, \quad  |\tilde f_\ell|\leq    g, \quad  |\langle \mu, f_\ell \rangle|\leq \langle \mu, g \rangle,\quad \Big|\sum_{j\in J}  \theta_j^n\,  \crr_j (f_\ell) \Big|\leq  2 g \quad\text{and}\quad |\crr_{j_0}(f_\ell)|\leq  c_J\,  g,
\end{equation}
where  for the  last inequality,  we  used that  the Vandermonde  matrix
$(\theta_j^n;   \,   j\in   J,    n\in   \{0,   \ldots,   |J|-1\})$   is
invertible.  Notice that  \reff{eq:unif-f-crit} holds  in particular  if
\reff{eq:hyp-crit} holds  for all $f\in  F$ and $\bF=(f_n, \,  n\in \N)$
takes   finitely  distinct   values  in   $F$  (\textit{i.e.}   the  set
$\{f_\ell; \ell\in \N\}\subset F$ is finite). The proof of the following
result is given in Section \ref{sec:proof-crit}.

\begin{theo}\label{theo:critical}
  Let  $X$  be  a  BMC   with  kernel  $\cp$  and  initial  distribution
  $\nu$. Assume that Assumptions  \ref{hyp:F} and \ref{hyp:F3} hold with
  $\alpha=1/\sqrt{2}$. We have the following convergence in distribution
  for  all  sequence  $\bF=(f_\ell,  \ell\in \N)$  of  elements  of  $F$
  satisfying    Assumptions    \ref{hyp:F3}     uniformly,    that    is
  \reff{eq:unif-f-crit}         for        some         $g\in        F$:
  \[ n^{-1/2} N_{n,  \emptyset}(\bF) \; \xrightarrow[n\rightarrow \infty
  ]{\text{(d)}}                          \;                         G,\]
  where  $G$ is  a Gaussian  real-valued random  variable with  variance
  $\scrit(\bF)$  given by  \reff{eq:scrit},  which is  well defined  and
  finite.
\end{theo}

\begin{rem}\label{rem:cr-F=ff}
For $f\in \cb(S)$, when it is well defined, we set:
\begin{equation}\label{eq:ST-crit}
\scrit_{\G}(f) = \sum_{j\in J} \langle \mu,
\cp(\crr_{j}(f) \sot \overline{\crr}_{j}(f)) \rangle
\quad\text{and}\quad
\scrit_{\TT}(f) = \scrit_{\G}(f) + 2\scrit_{\TT,2}(f),
\end{equation}
where
\[
  \scrit_{\TT,2}(f) 
=   \sum_{j\in J} \inv{\sqrt{2}\, \theta_j -1} \langle \mu, \cp(\crr_{j}(f) \sot
\overline{\crr}_{j}(f)) \rangle .
\]
 
If we take $\bF = (f,0,0, \ldots)$, we have $\scrit(\bF)=\scrit_{\G}(f)$. If we take $ \bff=(f,f, \ldots)$, the infinite sequence of the same function $f$, we have $\scrit(\bff)=2\scrit_{\TT}(f) $. 
\end{rem}

As a direct consequence of Remarks \ref{rem:cr-F=ff}  and
\ref{rem:simpleN0n}, and the more general Theorem
\ref{theo:critical}, we get
the following result. The proof which mimic the proof of Corollary
\ref{cor:subcritical} is left to the reader.

\begin{cor}\label{cor:critical}
Let   $X$  be  a  BMC   with  kernel  $\cp$  and  initial  distribution
$\nu$ such that Assumptions  \ref{hyp:F} and \ref{hyp:F3}
are in force with 
$\alpha= 1/\sqrt{2}$. 
Let  $f\in F$.
Then, we  have the  following convergence  in
  distribution:
\begin{equation*}
(n|\GG_{n}|)^{-1/2}M_{\GG_{n}}(\tilde{f}) \; \xrightarrow[n\rightarrow
\infty ]{\text{(d)}} \; G_{1},
\quad \text{and} \quad
  (n|\TT_{n}|)^{-1/2}M_{\TT_{n}}(\tilde{f}) \; \xrightarrow[n\rightarrow
  \infty ]{\text{(d)}} \; G_{2},  
\end{equation*}
where $G_{1}$ and $G_2$ are centered  Gaussian real-valued random
variables with respective variance
 $\scrit_{\GG}(f)$ and $\scrit_{\TT}(f)$
given in \reff{eq:ST-crit}, which are well defined and
finite.
\end{cor}

\begin{rem}
  We stress  that the variances $\scrit_{\TT}(f)$  and $\scrit_{\GG}(f)$
  can take the value 0. This is the case in particular if the projection
  of $f$ on the eigenspace corresponding to the eigenvalues $\alpha_{j}$
  equal 0  for all $j  \in J$. In the  symmetric BAR model  developed in
  Example \ref{ex:sBAR} where $J$ is reduced to a singleton and the
  projector is given by \eqref{eq:R-symBAR},  we deduce that if
  $a=\alpha=1/\sqrt{2}$  then  $\scrit_{\GG}(f)=\scrit_{\TT}(f)=0$ if 
  $\E[G\, f(\sigma_a  G)]=0$, where $G$ is  a standard $\cn(0,  1)$ Gaussian
  random variable. This is in particular  the case if $f$ is even.
\end{rem}


\subsection{The super-critical case $2\alpha^2>1$}
\label{sec:main-res-super}

We consider  the super-critical  case $\alpha\in (1/\sqrt{2}, 1)$. We shall  assume that Assumption \ref{hyp:F3}
holds. Recall \reff{eq:hyp-crit} with  the   eigenvalues  $\{\alpha_j=\theta_j \alpha, j\in J\} $ of $\cq$, with modulus equal to $\alpha$ (\textit{i.e.}
$|\theta_j|=1$) and the projector  $\crr_j$ on
the eigen-space  associated to eigenvalue  $\alpha_j$. 
Recall  that the sequence $(\beta_n,
n\in \N)$ in Assumption \ref{hyp:F3}
can (and will
) be chosen non-increasing and bounded from above by 1.

We shall consider the filtration $\ch=(\ch_n, n\in )$ defined by $\ch_{n} = \sigma(X_{i},i\in\TT_{n})$.  The next lemma, whose the proof is given in Section \ref{sec:martingaleSuC}, exhibits martingales related to the projector  $\crr_j$.

\begin{lem}\label{lem:martingale}
Let $X$ be a BMC with kernel $\cp$ and initial distribution
$\nu$. Assume that Assumption \ref{hyp:F} and \ref{hyp:F3} hold with
$\alpha\in (1/\sqrt{2}, 1)$ in
\reff{eq:hyp-crit}. Then, for all $j\in J$ and $f\in F$, the sequence
$M_{j}(f)=(M_{n,j}(f), n\in \N)$, with
\begin{equation*}
M_{n,j}(f) = (2 \alpha_j) ^{-n} \, M_{\G_n} (\crr_j(f)), 
\end{equation*}
is a $\ch$-martingale   which converges a.s. and in $L^{2}$ to a  random
variable, say  $M_{\infty,j}(f)$.
\end{lem}

Now, we state the main result of this section, whose proof is given in
Section \ref{sec:proof-super-thm}. Recall that
$\theta_j=\alpha_j/\alpha$ and $|\theta_j|=1$ and $M_{\infty , j}$ is
defined in Lemma \ref{lem:martingale}. 
\begin{theo}\label{theo:super-critical}
Let $X$ be a BMC with kernel $\cp$ and initial distribution
$\nu$. Assume that Assumptions \ref{hyp:F} and \ref{hyp:F3} hold with  $\alpha \in (1/\sqrt{2},
1) $ in  \reff{eq:hyp-crit}. 
We have the following convergence in probability for all sequence
$\bF=(f_\ell, \ell\in \N)$ of elements of 
$F$ satisfying Assumptions  \ref{hyp:F3} uniformly, that is
\reff{eq:unif-f-crit} holds for some $g\in F$:
\[
(2\alpha^{2})^{-n/2} N_{n, \emptyset}(\bF)
- \sum_{\ell\in \N} (2\alpha)^{-\ell} \sum_{j\in J} \theta_j^{n-\ell}
M_{\infty  ,j}(f_\ell)  
\; \xrightarrow[n\rightarrow \infty ]{\P}  \; 0. 
\]
\end{theo}
\begin{rem}\label{rem:theo-Scrit}
We stress that if for all $\ell \in \NN$, the orthogonal projection of $f_{\ell}$ on the eigen-spaces corresponding to the eigenvalues $1$ and $\alpha_{j}$, $j \in J$, equal 0, then $M_{\infty,j}(f_{\ell}) = 0$ for all $j\in J$ and in this case, we have
\begin{equation*}
(2\alpha^{2})^{-n/2} N_{n, \emptyset}(\bF) \; \xrightarrow[n\rightarrow \infty ]{\P}  \; 0. 
\end{equation*}
\end{rem}
 
As a direct consequence  of Theorem \ref{theo:super-critical} and Remark
\ref{rem:simpleN0n},  we  deduce  the  following  results.  Recall  that
$\tilde f= f -\langle \mu, f \rangle$.
\begin{cor}\label{cor:CtheoSC}
Under the assumptions of Theorem \ref{theo:super-critical}, we have for all $f \in F$:
\begin{align*}
(2\alpha)^{-n} M_{\TT_{n}}(\tilde{f}) - \sum_{j\in J} \theta_{j}^{n}(1 -
  (2\alpha\theta_{j})^{-1})^{-1}M_{\infty,j}(f) 
& \inP 0 \\
(2\alpha)^{-n} M_{\GG_{n}}(\tilde{f}) - \sum_{j\in J} \theta_{j}^{n}
  M_{\infty,j}(f) 
& \inP 0.
\end{align*}
\end{cor}
\begin{proof}
  We first  take $\bF = (f,f,\ldots)$  and next $\bF =  (f,0,\ldots)$ in
  Theorem \ref{theo:super-critical}, and then use \reff{eq:def-NOf}.
\end{proof}

We directly deduce the following two Corollaries.
\begin{cor}
   \label{cor:super-crit-d1}
Under the hypothesis of Theorem \ref{theo:super-critical}, if $\alpha$
is the only eigen-value of $\cq$ with modulus equal to $\alpha$ (and
thus $J$ is
reduced to a singleton), then we have:
\[
(2\alpha^{2})^{-n/2} N_{n, \emptyset}(\bF)
\; \xrightarrow[n\rightarrow \infty ]{\P}  \;
 \sum_{\ell\in \N} (2\alpha)^{-\ell}  M_{\infty }(f_\ell) ,
\]
where, for $f\in F$,  $M_{\infty }(f)=\lim_{n\rightarrow \infty } (2\alpha)^{-n}
M_{\G_n} (\crr (f))$, and $\crr$ is the projection on the eigen-space
associated to the eigen-value $\alpha$. 
\end{cor}

The next Corollary is a direct consequence of Corollary
\ref{cor:super-crit-d1}. 
\begin{cor}
   \label{cor:super-crit-d1-2}
Let $X$ be a BMC with kernel $\cp$ and initial distribution
$\nu$. Assume that Assumption \ref{hyp:F} and \ref{hyp:F3} hold with  $\alpha \in (1/\sqrt{2},
1) $ in  \reff{eq:hyp-crit}. Assume $\alpha$
is the only eigen-value of $\cq$ with modulus equal to $\alpha$ (and
thus $J$ is
reduced to a singleton), then we have for $f\in F$:
\begin{equation*}
  (2\alpha)^{-n}M_{\GG_{n}}(\tilde{f}) \inP M_{\infty}(f)
  \quad \text{and} \quad
  (2\alpha)^{-n}M_{\TT_{n}}(\tilde{f}) \inP \frac{2\alpha}{2\alpha - 1}M_{\infty}(f) ,
\end{equation*}
where $M_{\infty}(f)$ is a random variable defined in Corollary \ref{cor:super-crit-d1}.
\end{cor}

\section{Proof of Theorem \ref{theo:subcritical}}
\label{sec:proof-sub}
Let $(p_n, n\in \N)$ be a non-decreasing sequence of elements of $\N^*$
such that, for all $\lambda>0$:
\begin{equation}
   \label{eq:def-pn}
p_n< n, \quad \lim_{n\rightarrow \infty } p_n/n=1
\quad\text{and}\quad \lim_{n\rightarrow \infty } n-p_n - \lambda \log(n)=+\infty .
\end{equation}
When there is no ambiguity, we write $p$ for $p_n$. 
\medskip

Let $i,j\in \T$. We write $i\preccurlyeq  j$ if $j\in i\T$. We denote by
$i\wedge j$  the most recent  common ancestor of  $i$ and $j$,  which is
defined  as   the  only   $u\in  \T$   such  that   if  $v\in   \T$  and
$ v\preccurlyeq i$, $v \preccurlyeq j$  then $v \preccurlyeq u$. We also
define the lexicographic order $i\leq j$ if either $i \preccurlyeq j$ or
$v0  \preccurlyeq i$  and $v1  \preccurlyeq j$  for $v=i\wedge  j$.  Let
$X=(X_i, i\in  \T)$ be  a $BMC$  with kernel  $\cp$ and  initial measure
$\nu$. For $i\in \T$, we define the $\sigma$-field:
\begin{equation*}\label{eq:field-Fi}
\cf_{i}=\{X_u; u\in \T \text{ such that  $u\leq i$}\}.
\end{equation*}
By construction,  the $\sigma$-fields $(\cf_{i}; \, i\in \T)$ are nested
as $\cf_{i}\subset \cf_{j} $ for $i\leq  j$.
\medskip

We define for $n\in \N$, $i\in \G_{n-p_n}$ and $\bF\in F^\N$ the
martingale increments:
\begin{equation}
   \label{eq:def-DiF}
\Delta_{n,i}(\bF)= N_{n,i}(\bF) - \E\left[N_{n,i}(\bF) |\,
  \cf_i\right]
\quad\text{and}\quad \Delta_n(\bF)=\sum_{i\in \G_{n-p_n}} \Delta_{n,i}(\bF).
\end{equation}
Thanks to \reff{eq:def-NiF}, we have:
\[
\sum_{i\in \G_{n-p_n}} N_{n, i}(\bF) 
= |\G_n|^{-1/2} \sum_{\ell=0}^{p_n}  M_{\G_{n-\ell}} (\tilde f_\ell)
= |\G_n|^{-1/2} \sum_{k=n-p_n}^{n}  M_{\G_{k}} (\tilde f_{n-k}).
\]
Using the branching Markov property, and \eqref{eq:def-NiF}, we get for
$i\in \G_{n-p_n}$:
\[
\E\left[N_{n,i}(\bF) |\,
  \cf_i\right]
=\E\left[N_{n,i}(\bF) |\,
  X_i\right]
= |\G_n|^{-1/2} \sum_{\ell=0}^{p_n}
\E_{X_i}\left[M_{\G_{p_n-\ell}}(\tilde f_\ell)\right].
\] 
We deduce from \reff{eq:nof-D} with $k=n-p_n$ that:
\begin{equation}
   \label{eq:N=D+R}
N_{n, \emptyset}(\bF)
= \Delta_n(\bF) + R_0(n)+R_1(n),
\end{equation}
with
\begin{equation}
   \label{eq:reste01}
R_0(n)= |\G_n|^{-1/2}\, \sum_{k=0}^{n-p_n-1}
  M_{\G_k}(\tilde f_{n-k})
\quad\text{and}\quad
R_1(n)= \sum_{i\in \G_{n-p_n}}\E\left[N_{n,i}(\bF) |\,
  \cf_i\right].
\end{equation}
We have the following elementary lemma. 
\begin{lem}
   \label{lem:cvR0}
Under the assumptions of Theorem \ref{theo:subcritical}, we have the following convergence:
\[
\lim_{n\rightarrow \infty } \E[R_0(n)^2] =0.
\]
\end{lem}
\begin{proof}
For all $k\geq 1$, we have:
\begin{align*}
 \E_x[M_{\G_{k}} (\tilde f_{n-k})^2] 
& \leq  2^{k}\,  g_{1}(x) 
+ \sum_{\ell=0}^{k-1}  2^{k+\ell}\alpha^{2\ell}\, \, \cq^{k-\ell -1}
  \left(\cp ( g_{2}\otimes g_{2})\right)(x)\\
& \leq  2^{k} \,  g_{1}(x) 
+   2^{k} \sum_{\ell=0}^{k-1} (2\alpha^2)^{\ell}  g_{3}(x)\\
& \leq  2^{k} \,  g_4(x) ,
\end{align*}
with $g_1, g_2, g_3, g_4\in F$ and where we used \reff{eq:Q2}, \reff{eq:unif-f} twice and \reff{eq:erg-bd} twice (with $f$ and $g$ replaced
by $2(g^2+ \langle  \mu, g \rangle^2)$ and $g_{1}$, and with  $f$ and
$g$ replaced by $g$ and $g_{2}$)  for the first inequality, \reff{eq:erg-bd} (with
$f$ and $g$  replaced by $\cp ( g_{2}\otimes  g_{2})$ and $g_{3}$)
for the second, and that $2\alpha^2<1$ and $g_4=g_1+ (1-2\alpha^2)^{-1}
g_3$ for the last. As $g_4\in F\subset L^1(\nu)$, this
implies that $\E[M_{\G_{k}} (\tilde f_{n-k})^2] \leq  c^2 2^{k}$ for some
finite constant $c$ which does not depend on $n$ or $k$. We can take  $c$ large
enough, so that this upper bound  holds also for $k=0$ and all $n\in
\N$, thanks to \reff{eq:unif-f}. We deduce that:
\begin{equation}
   \label{eq:majo-R0}
\E[R_0(n)^2]^{1/2}
\leq  |\G_n|^{-1/2}\, \sum_{k=0}^{n-p-1} \E[M_{\G_{k}} (\tilde f_{n-k})^2]^{1/2}
\leq  c\, 2^{-n/2} \sum_{k=0}^{n-p-1} 2^{k/2}
\leq  3c\, 2^{- p/2}.
\end{equation}
Use that $\lim_{n\rightarrow \infty } p=\infty $ to conclude. 
\end{proof}

We have the following lemma. 
\begin{lem}
   \label{lem:cvR1}
Under the assumptions of Theorem \ref{theo:subcritical}, we have the following convergence:
\[
\lim_{n\rightarrow \infty } \E\left[R_1(n)^2\right]=0.
\]
\end{lem}

\begin{proof}
We  set  for $p\geq  \ell \geq 0$:
\begin{equation}
   \label{eq:def-R1ln}
R_1(\ell,n)=\sum_{i\in \G_{n-p}}\E\left[N_{n,i}^\ell(f_\ell) |\,
  \cf_i\right],
\end{equation}
so that, thanks to \eqref{eq:def-NiF},  $R_1(n)=\sum_{\ell=0}^{p}R_1(\ell,n)$. 
We have for $i\in \G_{n-p}$:
\begin{equation*}
   \label{eq:ENlni}
|\G_n|^{1/2} \E\left[N_{n,i}^\ell(f_\ell) |\,
  \cf_i\right] = \E\left[M_{i\G_{p-\ell}}(\tilde f_\ell)
  |X_i\right]=\E_{X_i} \left[M_{\G_{p-\ell}}(\tilde f_\ell)\right]=
 |\G_{p-\ell}|\, \cq^{p-\ell} \tilde f_\ell(X_i),
\end{equation*}
where we used 
definition \reff{eq:def-Nnil} of $N_{n,i}^\ell$ for the first equality,
the Markov property of $X$ for the second and  \reff{eq:Q1} for
the third. We deduce that:
\[
R_1(\ell,n)= |\G_n|^{-1/2} \,  |\G_{p-\ell}|\, M_{\G_{n-p}}
(\cq^{p-\ell} \tilde f_\ell).
\]
Using \reff{eq:Q2}, we get:
\begin{align*}
\E_x\left[R_1(\ell, n)^2\right] 
&= |\G_n|^{-1} \,  |\G_{p-\ell}|^{2}\, 
\E_x\left[\left(M_{\G_{n-p}} (\cq^{p-\ell} \tilde f_\ell)\right)^2\right]\\
&=  |\G_n|^{-1} \,  |\G_{p-\ell}|^{2}\, 2^{n-p} \,  \cq^{n-p}
\big(  (\cq^{p-\ell} \tilde f_\ell) ^2\big) (x) \\ 
&\hspace{1cm} +  |\G_n|^{-1} \,  |\G_{p-\ell}|^{2}\, \sum_{k=0}^{n-p-1}
  2^{n-p+k}\, \cq^{n-p-k-1} \left( \cp \left( \cq^{k+p-\ell} \tilde
  f_\ell \otimes^2 \right)\right)(x). 
\end{align*}
We deduce that:
\begin{align}
\nonumber
\E_x\left[R_1(\ell, n)^2\right] 
&\leq  \alpha^{2(p-\ell)} 2^{p-2\ell} \,  \cq^{n-p} (g^2)(x) 
 + 2^{p-2\ell} \, \sum_{k=0}^{n-p-1}
  \alpha^{2(k+p-\ell)}\, 2^k \, 
\cq^{n-p-k-1} \left( \cp \left( g \otimes g \right)\right) \\
\label{eq:R1-ineg-sub}
&\leq \alpha^{2(p-\ell)} 2^{p-2\ell} \, \left(g_{1}(x)+ \sum_{k=0}^{n-p-1} (2
  \alpha^2)^k g_{2}(x)\right) \\
\nonumber
&\leq   (2 \alpha^2)^p\, 
(2\alpha)^{-2\ell} \, g_{3}(x),
\end{align}
with  $g_1, g_2, g_3\in  F$ and  where  we  used  \reff{eq:unif-f} for  the  first
inequality, \reff{eq:erg-bd} twice  (with $f$ and $g$  replaced by $g^2$
and $g_{1}$ and by $\cp \left( g \otimes g \right)$ and $g_{2}$) for the
second,    and    that    $2\alpha^2<1$     for    the    last.    Since
$g_3\in       F\subset       L^1(\nu)$,      this       gives       that
$  \E\left[R_1(\ell, n)^2\right]\leq  (2 \alpha^2)^p  (2\alpha)^{-2\ell}
\langle \nu, g_{3} \rangle$.  We deduce that:
\[
\E\left[R_1(n)^2\right]^{1/2}
\leq  \sum_{\ell=0}^{p} \E\left[R_1(\ell, n)^2\right]^{1/2}
\leq  a_{1, n}\, \langle \nu, g_3 \rangle^{1/2},
\]
with the sequence $(a_{1, n},n\in \N)$ defined by:
\begin{equation*}
   \label{eq:def-a1n}
a_{1, n}= (2\alpha^2)^{p/2} \sum_{\ell=0}^{p}
(2\alpha)^{-\ell}. 
\end{equation*}
Notice the sequence $(a_{1, n},n\in \N)$ converges to 0 since $\lim_{n\rightarrow\infty } p
=\infty $, $2\alpha^2<1$  and 
\begin{equation*}\label{eq:sum2alpha-}
\sum_{\ell=0}^{p} (2\alpha)^{-\ell}\leq  
\begin{cases} 
2\alpha/(2\alpha -1) & \text{if $2\alpha>1$},\\ 
p+1 & \text{if $2\alpha=1$},\\ 
(2\alpha)^{-p}/(1-2\alpha) & \text{if $2\alpha<1$}. 
\end{cases}
\end{equation*}
We conclude that $\lim_{n\rightarrow \infty }
\E\left[R_1(n)^2\right]=0$.
 \end{proof}

We now study the bracket of $\Delta_n$:
\begin{equation}
   \label{eq:def-Vn-subc}
V(n)= \sum_{i\in \G_{n-p_n}} \E\left[ \Delta_{n, i}
  (\bF)^2|\cf_i\right]. 
\end{equation}
Using \reff{eq:def-NiF} and \reff{eq:def-DiF}, we write:
\begin{equation}
   \label{eq:def-V}
V(n) = |\G_n|^{-1} \sum_{i\in \G_{n-p_n}} \E_{X_i}\left[
    \left(\sum_{\ell=0}^{p_n} M_{\G_{p_n-\ell}}(\tilde f_\ell) \right)^2
  \right]-R_2(n)=V_1(n) +2V_2(n) - R_2( n),
\end{equation}
with:
\begin{align*}
V_1(n)
& =   |\G_n|^{-1} \sum_{i\in \G_{n-p_n}} 
\sum_{\ell=0}^{p_n} \E_{X_i}\left[
   M_{\G_{p_n-\ell}}(\tilde f_\ell) ^2  \right] ,\\
V_2(n)
& =  |\G_n|^{-1} \sum_{i\in \G_{n-p_n}} \sum_{0\leq \ell<k\leq p_n
  } \E_{X_i}\left[ 
   M_{\G_{p_n-\ell}}(\tilde f_\ell)  M_{\G_{p_n-k}}(\tilde f_k) 
  \right], \\
R_2( n)
&=\sum_{i\in \G_{n-p_n}} \E\left[ N_{n,i} (\bF)
  |X_i \right] ^2.
\end{align*}

\begin{lem}
   \label{lem:cvR2}
Under the assumptions of Theorem \ref{theo:subcritical}, we have the following convergence:
\[
\lim_{n\rightarrow \infty } \E\left[R_2(n)\right]=0.
\]
\end{lem}

\begin{proof}
We define the sequence $(a_{2, n}, n\in \N)$ for $n\in \N$ by: 
\begin{equation*}
   \label{eq:def-a2n}
a_{2, n}= 2^{-p} 
 \left(\sum_{\ell=0}^{p}  (2\alpha)^{\ell} \right)^2.
\end{equation*}
Notice that the sequence $(a_{2, n}, n\in \N)$   converges to 0 since 
$\lim_{n\rightarrow \infty }p=\infty $, $2\alpha^2<1$  and 
\begin{equation*}
   \label{eq:sum2a}
\sum_{\ell=0}^{p}
(2\alpha)^{\ell}\leq  
\begin{cases}
(2\alpha)^{p+1}/(2\alpha-1) & \text{if $2\alpha>1$},\\
p+1 & \text{if $2\alpha=1$},\\
   1/(1-2\alpha ) & \text{if $2\alpha<1$}.
\end{cases}
\end{equation*}

 We now compute $ \E_x\left[R_2(n)\right]$. 
\begin{align}
\nonumber
  \E_x\left[R_2(n)\right]
 &= |\G_n|^{-1} \, \sum_{i\in \G_{n-p}} \E_x\left[
\E_x\left[\sum_{\ell=0}^{p} M_{i\G_{p-\ell}}(\tilde f_\ell)
  |X_i\right]^2
\right]\\
\nonumber
 &= |\G_n|^{-1} \, \sum_{i\in \G_{n-p}} \E_x\left[
\left(\sum_{\ell=0}^{p} \E_{X_i} \left[M_{\G_{p-\ell}}(\tilde f_\ell)
 \right]\right)^2
\right]\\
\nonumber
&= |\G_n|^{-1} \,  |\G_{n-p}|\, 
\cq^{n-p} \left(\Big(\sum_{\ell=0}^{p}  |\G_{p-\ell}|\, \cq^{p-\ell} \tilde
  f_\ell\Big)^2\right)(x)
\\
\label{eq:majoR2-sub-ineg}
&\leq 2^{-p} 
 \left(\sum_{\ell=0}^{p}  (2\alpha)^{p-\ell} \right)^2\cq^{n-p}(g^2) (x)\\
\nonumber
&\leq a_{2, n}\,  g_1 (x),
\end{align}
with $g_1\in F$ and where we used the definition of $N_{n,i}(\bF)$ for the first equality,
the Markov property of $X$ for the second,  
\reff{eq:Q1} for the
third,  \reff{eq:unif-f} for the first inequality,
 and \reff{eq:erg-bd} (with $f$ and $g$
replaced by $g^2$ and $g_{1}$) for the last.
We conclude that $\lim_{n\rightarrow \infty }
\E\left[R_2(n)\right]=0$, 
using that $\langle \nu, g_1 \rangle$ if finite as $g_1\in F\subset
L^1(\nu)$. 
\end{proof}

We have the following technical lemma.
\begin{lem}
   \label{lem:cvV2}
Under the assumptions of Theorem \ref{theo:subcritical}, we have that
$\ssub_2(\bF)$ defined in \reff{eq:S2} is well defined and finite, and
that a.s. 
$ \lim_{n\rightarrow \infty } V_2(n) =\ssub_2(\bF)<+\infty $. 
\end{lem}

\begin{proof}
   Using \reff{eq:Q2-bis}, we get:
\begin{equation}
   \label{eq:decom-V2}
V_2(n)= V_5(n)+ V_6(n),
\end{equation}
with
\begin{align*}
    V_5(n)
&=  |\G_n|^{-1} \sum_{i\in \G_{n-p}} 
\sum_{0\leq \ell<k\leq  p } 2^{p-\ell} \cq^{p-k} \left( \tilde f_k
  \cq^{k-\ell} \tilde f_\ell\right)(X_i),\\
V_6(n)
&=     |\G_n|^{-1} \sum_{i\in \G_{n-p}} 
\sum_{0\leq \ell<k<  p }
\sum_{r=0}^{p-k-1}  2^{p-\ell+r} \,
  \cq^{p-1-(r+k)}\left(\cp\left(\cq^r \tilde f_k \sot
\cq  ^{k-\ell+r} \tilde f_\ell \right)\right)(X_i).
\end{align*}
We consider the term $V_6(n)$. We have:
\begin{equation}
   \label{eq:V6}
V_6(n)=|\G_{n-p}|^{-1} M_{\G_{n-p}} (H_{6,n}),
\end{equation}
with:
\begin{equation*}
   \label{eq:def-Hn-V6}
H_{6,n}=\sum_{\substack{0\leq \ell< k \\ r\geq 0}}  h_{k,\ell,r}^{(n)}\,
\ind_{\{r+k<  p\}} 
\,\text{ and }\,
h_{k, \ell,r}^{(n)} =  2^{r-\ell} \,
  \cq^{p-1-(r+k)}\left(\cp\left(\cq^r \tilde f_k \sot
\cq  ^{k-\ell+r} \tilde f_\ell \right)\right).
\end{equation*}
Using              \reff{eq:geom-erg}             and              since
$\cp(\cq^r(F)\otimes       \cq^{k-\ell+r}(F))\subset        F$       and
$\lim_{n\rightarrow \infty } p=+\infty $, we have that:
\[
\lim_{n\rightarrow \infty } h_{k,\ell,r}^{(n)} = h_{k,\ell,r},
\]
where the constant $h_{k,\ell,r}$ is equal to $2^{r-\ell} \langle \mu,
\cp\left(\cq^r \tilde 
    f_k \sot \cq^{k-\ell+r} \tilde f_\ell\right)\rangle$. Using \eqref{eq:geom-erg}, we also have that:
\begin{align*}
   |h_{k, \ell,r}^{(n)}|
&\leq   2^{r-\ell}\, \alpha^{k-\ell+2r}\,  \cq^{p-1-(r+k)}\left(\cp\left(g \otimes
  g\right)\right)\\
&\leq   2^{r-\ell}\, \alpha^{k-\ell+2r}\,  g_*,
\end{align*}
with $g_*\in F$ (which does not depend on $n, r, k$ and $\ell$) and where we used \reff{eq:unif-f} for the first inequality and
\reff{eq:erg-bd} (with $f$ and $g$ replaced by $\cp\left(g \otimes
  g\right)$ and $g_*$). Taking the limit, we also deduce that:
\[
    |h_{k, \ell,r}|\leq    2^{r-\ell}\, \alpha^{k-\ell+2r}\,  g_*.
\]
Define the constant 
\begin{equation*}
   \label{eq:def-H6}
H_6(\bF)=\sum_{\substack{0\leq \ell< k \\ r\geq 0}} h_{k, \ell,r}
=\sum_{\substack{0\leq \ell< k \\ r\geq 0}} 
2^{r-\ell} \langle \mu,
\cp\left(\cq^r \tilde 
    f_k \sot \cq^{k-\ell+r} \tilde f_\ell\right)\rangle
\end{equation*}
 which is finite as 
\begin{equation}
   \label{eq:sum2alpha-b}
\sum_{0\leq \ell< k, \, r\geq 0}
2^{r-\ell}\, \alpha^{k-\ell+2r}=\frac{2\alpha}{(1-\alpha)(1-2\alpha^2)}<+\infty
.
\end{equation}
  Using \reff{eq:geom-erg} (with $f$ and $g$ replaced by
$\cp\left(\cq^r \tilde
    f_k \sot \cq^{k-\ell+r} \tilde f_\ell\right)$ and $g_{k,\ell,r}$), we deduce that:
\[
| h_{k, \ell,r}^{(n)} - h_{k, \ell,r}|\leq  2^{r-\ell} \alpha^{p-1-(r+k)} g_{k,\ell,r}.
\]
Set $r_0\in \N^*$  and $g_{r_0}=\sum_{0\leq
  \ell< k;\, r\geq 0;\, k\vee r\leq r_0} g_{k, \ell,r}$. Notice that $g_{r_0}$
belongs to $F$ and is non-negative. Furthermore, we have:
\begin{align*}
|H_{6,n} -H_6(\bF)|
&\leq  \sum_{\substack{0\leq \ell< k\\ r\geq 0\\  k\vee r\leq r_0}}  2^{r-\ell}
\alpha^{p-1-(r+k)} g_{r_0}
+ \sum_{\substack{0\leq \ell< k\\ r\geq 0\\ r \vee k>r_0}}\left( |h_{k,
  \ell,r}^{(n)}|    \,  \ind_{\{r+k<  p\}} 
+ |h_{k,\ell,r}|\right)\\
&\leq (r_0+1)^2\, 2^{r_0+1} \, \alpha^{p-1-2r_0} \,  g_{r_0} +
\gamma_1(r_0)g_*,
\end{align*}
with
\[
\gamma_1(r_0)= \sum_{\substack{0\leq \ell< k\\ r\geq 0\\ r  \vee 
   k>r_0}} 2^{r-\ell} \, \alpha^{ k-\ell+2r}.
\]
Using \reff{eq:lfgn-G} with $n$ replaced by $n-p$ and $f$ replaced by $g_*$ and $g_{r_0}$, and
 that $\lim_{n\rightarrow \infty } \alpha ^p=0$ as well as
 $\lim_{n\rightarrow \infty } n-p=\infty $, we deduce
that:
\[
\limsup_{n\rightarrow \infty } 
|\G_{n-p}|^{-1}\, M_{\G_{n-p}}(|H_{6,n}-H_6(\bF)|)\leq  
\gamma(r_0) 
\langle \mu, g_* \rangle.
\]
Thanks to \reff{eq:sum2alpha-b}, we get by dominated convergence that
$\lim_{r_0 \rightarrow \infty  } \gamma_1(r_0)=0$. This implies  that:
\[
\lim_{n\rightarrow \infty } 
|\G_{n-p}|^{-1}\, M_{\G_{n-p}}(|H_{6,n}-H_6(\bF)|)=0.
\]
Since $|\G_{n-p}|^{-1}\, M_{\G_{n-p}}(\cdot)$ is a probability measure, we
deduce from \reff{eq:V6} that a.s.:
\[
\lim_{n\rightarrow\infty } V_6(n)=\lim_{n\rightarrow\infty }
|\G_{n-p}|^{-1}\, M_{\G_{n-p}}(H_{6,n}) 
=H_6(\bF)=\sum_{\substack{0\leq \ell< k \\ r\geq 0} } 2^{r-\ell}
\langle \mu, \cp\left(\cq^r \tilde 
    f_k \sot \cq^{k-\ell+r} \tilde f_\ell\right)\rangle.
\]
Similarly, we get that  a.s. $\lim_{n\rightarrow\infty } V_5(n) = H_5(\bF)$, 
with the finite constant $H_5(\bF)$ defined by:
\begin{equation*}
   \label{eq:def-H5}
H_5(\bF)=
\sum_{0\leq \ell< k}
2^{-\ell} \langle \mu,   \tilde f_k \cq^{k-\ell} \tilde f_\ell\rangle.
\end{equation*}
Notice that $ \ssub_2(\bF)=H_5(\bF)+H_6(\bF)$ is finite thanks to
\reff{eq:unif-f} and 
\reff{eq:sum2alpha-b}. This finishes the proof. 
\end{proof}

Using similar arguments as in the proof of Lemma \ref{lem:cvV2}, we get
the following result. 
\begin{lem}
   \label{lem:cvV1}
Under the assumptions of Theorem \ref{theo:subcritical}, we have that $
\ssub_1(\bF)$ in \reff{eq:S1} is well defined and finite, and that 
a.s. $\lim_{n\rightarrow \infty } V_1(n) =\ssub_1(\bF)$.
\end{lem}

\begin{proof}
Using \reff{eq:Q2}, we get:
\begin{equation}\label{eq:DV4nsub}
V_1(n)= V_3(n)+ V_4(n),
\end{equation}
with
\begin{align*}
V_3(n) &=  |\G_n|^{-1} \sum_{i\in \G_{n-p}} \sum_{\ell=0}^p 2^{p-\ell}\,
         \cq^{p-\ell} (\tilde f_\ell^2)(X_i),\\ 
V_4(n) &=     |\G_n|^{-1} \sum_{i\in \G_{n-p}} \sum_{\ell=0}^{p-1}\,
         \sum_{k=0}^{p-\ell -1} 2^{p-\ell+k}
         \,\cq^{p-1-(\ell+k)}\left(\cp\left(\cq^k \tilde f_\ell
         \otimes^2\right)\right)(X_i).  
\end{align*}
We consider the term $V_4(n)$. We have:
\begin{equation}
   \label{eq:V4}
V_4(n)=|\G_{n-p}|^{-1} M_{\G_{n-p}} (H_{4,n}),
\end{equation}
with:
\begin{equation}
   \label{eq:def-H4n-sub-pt}
H_{4,n}=\sum_{\ell\geq 0, \, k\geq 0} h_{\ell,
  k}^{(n)}\,  \ind_{\{\ell+k<  p\}}
\quad\text{and}\quad
h_{\ell,
  k}^{(n)} = 2^{k-\ell}\,  \cq^{p-1-(\ell+k)} \left(\cp\left(\cq^k \tilde
    f_\ell \otimes^2 \right)\right).
\end{equation}
Using \reff{eq:geom-erg}, we have that:
\[
\lim_{n\rightarrow \infty } h_{\ell,
  k}^{(n)} = h_{\ell,
  k},
\]
where the constant $h_{\ell,
  k}$ is equal to $2^{k-\ell} \langle \mu, \cp\left(\cq^k \tilde
    f_\ell \otimes^2\right)\rangle$. We also have that:
\begin{align*}
   |h_{\ell,
  k}^{(n)}|
&\leq   2^{k-\ell}\, \alpha^{2k}\,  \cq^{p-1-(\ell+k)}\left(\cp\left(g \otimes
  g\right)\right)\\
&\leq   2^{k-\ell}\, \alpha^{2k}\,  g_*,
\end{align*}
with $g_*\in F$ (which  does not depend on $n, \ell$  and $k$) and where
we used  \reff{eq:unif-f} for the first  inequality and \reff{eq:erg-bd}
(with  $f$  and  $g$  replaced  by  $\cp\left(g  \otimes  g\right)$  and
$g_*$). Taking the limit, we also deduce that:
\[
    |h_{\ell,
  k}|\leq   2^{k-\ell}\, \alpha^{2k}\,  g_*.
\]
Define the constant 
\begin{equation*}
   \label{eq:def-H4}
H_4(\bF)=\sum_{\ell\geq 0, \, k\geq 0} h_{\ell,k},
\end{equation*}
which
is finite as:
\begin{equation}
   \label{eq:sum2alpha2}
\sum_{\ell\geq 0,\, k\geq 0} 2^{k-\ell}\alpha^{2k}=2/(1-2\alpha^2)<+\infty .
 \end{equation}   
Using \reff{eq:geom-erg} (with $f$ and $g$ replaced by $\cp\left(\cq^k \tilde f_\ell \otimes^2\right)$ and $g_{\ell,k}$), we deduce that:
\[
| h_{\ell,
  k}^{(n)} - h_{\ell,
  k}|\leq  2^{k-\ell} \alpha^{p-1-(\ell+k)} g_{\ell,k},
\]
Set $r_0\in \N$ and $g_{r_0}=\sum_{\ell\vee  k\leq r_0} g_{\ell,
  k}$. Notice that $g_{r_0}$
belongs to $F$. Furthermore, we have:
\begin{align*}
|H_{4,n} -H_4(\bF)|
&\leq  \sum_{\ell\vee k\leq r_0}  2^{k-\ell}
\alpha^{p-1-(\ell+k)} g_{r_0}
+ \sum_{\ell \vee k>r_0} \left(|h_{\ell, k}^{(n)}|    \,  \ind_{\{\ell+k\leq
  p-1\}}
+  |h_{\ell, k}|\right)\\
&\leq (r_0+1)^2 2^{r_0} \, \alpha^{p-1-2r_0} \,  g_{r_0} + \gamma_2(r_0) g_*,
\end{align*}
with $\gamma_2(r_0) =2\sum_{\ell \vee
  k>r_0} 2^{k-\ell} \, \alpha^{2 k}$. 
Using \reff{eq:lfgn-G}  with $n$ replaced by $n-p$ and  $f$ replaced by $g_*$ and $g_{r_0}$,  and
 that $\lim_{n\rightarrow \infty } \alpha ^p=0$ as well as
 $\lim_{n\rightarrow \infty } n-p=\infty $, we deduce
that:
\[
\limsup_{n\rightarrow \infty } 
|\G_{n-p}|^{-1}\, M_{\G_{n-p}}(|H_{4,n}-H_4(\bF)|)\leq  
\gamma_2(r_0) \langle \mu, g_* \rangle.
\]
Thanks to \reff{eq:sum2alpha2}, we get by dominated convergence that
$\lim_{r_0 \rightarrow \infty  } \gamma_2(r_0)=0$. 
We deduce that:
\[
\lim_{n\rightarrow \infty } 
|\G_{n-p}|^{-1}\, M_{\G_{n-p}}(|H_{4,n}-H_4(\bF)|)=0.
\]
Since $|\G_{n-p}|^{-1}\, M_{\G_{n-p}}(\cdot)$ is a probability measure, we
deduce from \reff{eq:V4} that a.s.:
\[
\lim_{n\rightarrow\infty } V_4(n)=\lim_{n\rightarrow\infty } |\G_{n-p}|^{-1}\, M_{\G_{n-p}}(H_{4,n})
=H_4(\bF)=\sum_{\ell\geq 0, \, k\geq 0} 
2^{k-\ell} \langle \mu, \cp\left(\cq^k \tilde
    f_\ell \otimes^2\right)\rangle.
\]
Similarly, we get that a.s. $\lim_{n\rightarrow\infty } V_3(n)=H_3(\bF)$
with the finite constant $H_3(\bF)$ defined by
\begin{equation*}
   \label{eq:def-H3}
H_3(\bF)=\sum_{\ell\geq 0} 
2^{-\ell} \langle \mu,   \tilde f_\ell^ 2\rangle.
\end{equation*}
Notice that $ \ssub_1(\bF)=H_3(\bF)+H_4(\bF)$ is finite thanks to
\reff{eq:unif-f} and 
\reff{eq:sum2alpha2}. This finishes the proof. 
\end{proof}

 The  next Lemma is a direct consequence of \reff{eq:def-V} and
 Lemmas \ref{lem:cvR2}, \ref{lem:cvV2} and \ref{lem:cvV1}.
\begin{lem}
   \label{lem:cvV}
Under the assumptions of Theorem \ref{theo:subcritical}, we have the
following convergence in probability $\lim_{n\rightarrow\infty }
V(n)=\ssub(\bF)$, where, with $\ssub_1(\bF)$ and  $\ssub_2(\bF)$ 
 defined by \reff{eq:S1} and \reff{eq:S2}:
\[
\ssub(\bF)=\ssub(\bF)=\ssub_1(\bF)+ 2\ssub_2(\bF). 
\]
\end{lem}

We now check the Lindeberg condition using a fourth moment condition. We set:
\begin{equation}
   \label{eq:def-R3}
R_3(n)=\sum_{i\in \G_{n-p_n}} \E\left[\Delta_{n,i}(\bF)^4\right].
\end{equation}

\begin{lem}
   \label{lem:cvG}
Under the assumptions of Theorem \ref{theo:subcritical}, we have that
$\lim_{n\rightarrow\infty } R_3(n)=0$.
\end{lem}
\begin{proof}
We have:
\begin{align}
\nonumber
 R_3(n)
&\leq   16 \sum_{i\in \G_{n-p}}
 \E\left[N_{n,i}(\bF)^4\right]\\
\nonumber
&\leq   16 (p+1)^3 \sum_{\ell=0}^p \sum_{i\in \G_{n-p}}
\E\left[N_{n,i}^\ell(\tilde f_\ell)^4\right],
\end{align}
where we used that $(\sum_{k=0}^r a_k)^4 \leq  (r+1)^3 \sum_{k=0}^r
a_k^4$ for the two inequalities (resp. with $r=1$ and $r=p$) and also
Jensen inequality and \reff{eq:def-DiF} for the first and 
\reff{eq:def-NiF} for the last. 
Using \reff{eq:def-Nnil}, we get:
\begin{equation*}
   \label{eq:def-R3hnl}
\E\left[N_{n,i}^\ell(\tilde f_\ell)^4\right]= |\G_n|^{-2} \E\left[h_{n,\ell}(X_i)\right], 
\quad\text{with}\quad
h_{n,\ell}(x)=\E_x\left[M_{\G_{p-\ell}}(\tilde f_\ell) ^4\right]. 
\end{equation*}
Thanks to the fourth moment bound given in Lemma \ref{lem:M4}, the
uniform bounds from \eqref{eq:unif-f} and the structural assumption
\ref{hyp:F}, it is easy to get there  exists $g_1\in F$ such  that for all
$n\geq p \geq \ell\geq 0$:
\begin{equation}
   \label{eq:majo-hnl-g1}
|h_{n,\ell}|\leq  2^{2(p-\ell)} g_1.
\end{equation}
We deduce that:
\begin{align*}
R_3(n) 
&\leq   16 n^3 \sum_{\ell=0}^p \sum_{i\in \G_{n-p}}
|\G_n|^{-2}\,  2^{2(p-\ell)} \E\left[g_1(X_i)\right]\\
&\leq   16 n^3 \,2^{-2(n-p)} \, \E\left[M_{\G_{n-p}}(g_1)\right]\\
&\leq   16 n^3 \, 2^{-(n-p)} \, \langle \nu, \cq^{n-p} g_1 \rangle,
\end{align*}
where we used \reff{eq:Q1} for the third inequality. Since $g_1$ belongs
to $F$, we deduce from \reff{eq:erg-bd} that $\cq^{n-p} g_1\leq  g_2$
for some $g_2\in F$ and all $n\geq  p\geq 0$. This gives that:
\[
R_3(n) \leq    16 n^3\, 2^{-(n-p)} \, \langle \nu, g_2 \rangle.
\]
This  ends the proof as $\lim_{n\rightarrow \infty } p=\infty $
and $\lim_{n\rightarrow\infty } n-p - \lambda \log(n)=+\infty $ for all $\lambda>0$.
\end{proof}

We can  now use  Theorem 3.2 and  Corollary 3.1, p.~58, and  the remark
p.~59  from  \cite{hh:ml}  to  deduce  from  Lemmas  \ref{lem:cvV}  and
\ref{lem:cvG} that  $\Delta_n(\bF)$ converges in distribution  towards a
Gaussian  real-valued   random  variable  with   deterministic  variance
$\ssub(\bF)$ given  by \reff{eq:ssub}.  Using  \reff{eq:N=D+R} and
Lemmas  \ref{lem:cvR0}  and  \ref{lem:cvR1},   we  then  deduce  Theorem
\ref{theo:subcritical}.

\section{Proof of Theorem \ref{theo:critical}}\label{sec:proof-crit}
We keep notation from Section \ref{sec:proof-sub}. Let $(p_n, n\in \N)$ be an increasing sequence of elements of $\N$ such that \reff{eq:def-pn} holds. When there is no ambiguity, we write $p$ for $p_n$. Recall the definitions of $\Delta_n(\bF)$ and $N_{n, \emptyset}(\bF)$ from \reff{eq:def-DiF} and \reff{eq:N=D+R}, as well as $R_0(n)$ and $R_1(n) $ from \reff{eq:reste01}. We have the following elementary lemma. 
\begin{lem}\label{lem:cvR0-crit}
Under the assumptions of Theorem \ref{theo:critical}, we have the following convergence: \[\lim_{n\rightarrow \infty } n^{-1} \E[R_0(n)^2] =0.\]
\end{lem}
\begin{proof}
Following the proof of Lemma \ref{lem:cvR0}, and using that
$2\alpha^2=1$ so that $\sum_{\ell=0}^{k-1} (2\alpha^2)^{\ell} =k$, we get
there exists some finite constant $c$ depending on $\bF$ such that  $\E[M_{\G_k} (\tilde f_{n-k})^2] \leq  c^2\,  (k+1) 2^k$ for all $k\geq
0$. This implies that:
\[
\E[R_0(n)^2]^{1/2}
\leq  |\G_n|^{-1/2}\, \sum_{k=0}^{n-p-1} \E[M_{\G_{k}} (\tilde f_{n-k})^2]^{1/2}
\leq  c\,  2^{-n/2} \sum_{k=0}^{n-p-1} \sqrt{k+1} \, 2^{k/2}
\leq   Cc\, \sqrt{n}\, 2^{- p/2}.
\]
Then use that $\lim_{n\rightarrow \infty } p/n=1 $ to conclude. 
\end{proof}

We have the following lemma. 
\begin{lem}
   \label{lem:cvR1-crit}
Under the assumptions of Theorem \ref{theo:critical}, we have the following convergence:
\[
\lim_{n\rightarrow \infty } n^{-1} \E\left[R_1(n)^2\right]=0.
\]
\end{lem}

\begin{proof}
Following the proof of Lemma \ref{lem:cvR1} with the same notations, and using that
$2\alpha^2=1$ so that $\sum_{k=0}^{n-p-1} (2\alpha^2)^{k} =n-p$
in \reff{eq:R1-ineg-sub}, we get
that there exists $g_3\in F$ such that 
 $ \E\left[R_1(\ell, n)^2\right]\leq (n-p+1)
(2\alpha)^{-2\ell} \langle \nu, g_3 \rangle$, where $R_1(\ell, n)$ is
defined in \reff{eq:def-R1ln}. 
As $2\alpha=\sqrt{2}$ and $R_1(n)=\sum_{\ell=0}^p R_1(\ell, n)$, we deduce that:
\[
\E\left[R_1(n)^2\right]^{1/2}
\leq  \sum_{\ell=0}^{p} \E\left[R_1(\ell, n)^2\right]^{1/2}
\leq 4 \sqrt{n-p+1} \, \langle \nu, g_3 \rangle^{1/2}  .
\]
Use that $\lim_{n\rightarrow\infty } p/n=1$ to conclude. 
\end{proof}

Recall $\Delta_{n} (\bF)$ defined  in \reff{eq:def-DiF}, and its bracket
defined                                                               by
$V(n)=     \sum_{i\in     \G_{n-p_n}}     \E\left[     \Delta_{n,     i}
  (\bF)^2|\cf_i\right]$ defined in \reff{eq:def-Vn-subc}. Recall, see
\reff{eq:def-V}, that $V(n)=V_1(n) +2V_2(n) - R_2( n)$. We study the
convergence of each term of the latter right hand side. 

\begin{lem}
   \label{lem:cvR2-crit}
Under the assumptions of Theorem \ref{theo:critical}, we have the following convergence:
\[
\lim_{n\rightarrow \infty } n^{-1/2} \E\left[R_2(n)\right]=0.
\]
\end{lem}

\begin{proof}
  Following   the  proof   of  Lemma   \ref{lem:cvR2}  with the same
  notations and   using  that
  $2\alpha^2=1$ so that $\sum_{\ell=0}^p  (2\alpha)^\ell \leq C 2^{p/2}$
  in        \reff{eq:majoR2-sub-ineg},         we        get        that
  $\E[R_2(n)]\leq C\, \langle  \nu, g_1 \rangle$, with  $g_1\in F$. This
  gives the result.
\end{proof}

Recall  $f_{k, \ell}^*$ defined in \reff{eq:def-fkl*}. 
 For $k, \ell,r\in \N$, we will  consider the
$\C$-valued functions  on $S^2$: 
\begin{equation}\label{eq:def-fkl}
f_{k,\ell,r} = \Big( \sum_{j\in J}
  \theta_j^{r}\,  \crr_j (f_k) \Big)
\sot \Big(\sum_{j\in J}  \theta_j^{r+k-\ell}\,  \crr_j (f_\ell) \Big)
\quad\text{and}\quad
f_{k, \ell,r}^\circ=f_{k, \ell,r} -  f_{k, \ell}^* .
\end{equation}

\begin{lem}\label{lem:cvV2-crit}
Under the assumptions of Theorem \ref{theo:critical}, we have that
a.s. $\lim_{n\rightarrow \infty } n^{-1} V_2(n) = \scrit _2(\bF)$ with
$\scrit _2(\bF)$ defined by \reff{eq:S2-crit} which is well defined and
finite. 
\end{lem}

\begin{proof}
We keep the decomposition \reff{eq:decom-V2} of $V_2(n)=V_5(n)+V_6(n)$ given in the proof of Lemma \ref{lem:cvV2}. We first consider the term $V_6(n)$ given in \reff{eq:V6} by:
\begin{equation}\label{eq:def-V6-crit}
 V_6(n)= |\G_{n-p}|^{-1} M_{\G_{n-p}} (H_{6,n}) ,
\end{equation}
with:
\begin{equation*}\label{eq:def-V6-crit-detail}
H_{6,n}=\!\!\sum_{\substack{0\leq \ell<k\leq  p\\ r\geq 0}}\!\!
h_{k,\ell,r}^{(n)}\,  \ind_{\{r+k<  p\}} 
\quad\text{and}\quad 
h_{k,\ell,r}^{(n)} =  2^{r-\ell} \, \cq^{p-1-(r+k)}\left(\cp\left(\cq^r
    \tilde f_k \sot \cq  ^{k-\ell+r} \tilde f_\ell \right)\right). 
\end{equation*}
We set 
\begin{equation*}\label{eq:Hntilde-crit}
\bar{H}_{6,n} = \sum_{0\leq \ell<k\leq  p; \, r\geq 0}  \bar h_{k,\ell,r}^{(n)}\,  \ind_{\{r+k<  p\}} 
\end{equation*}
where  for $0\leq \ell<k\leq p$ and $0\leq r<p-k$: 
\[
\bar  h_{k,\ell,r}^{(n)} =  2^{r-\ell} \,  \alpha^{k-\ell +2r} \,
\cq^{p-1-(r+k)}(\cp f_{k, \ell,r})
=  2^{-(k+\ell)/2}\, \cq^{p-1-(r+k)}(\cp f_{k, \ell,r}), 
\]
where  we used  that  $2\alpha^2=1$. We have:
\begin{align*}
|h_{k, \ell,r}^{(n)} - \bar  h_{k, \ell,r}^{(n)}| 
&\leq  2^{r-\ell} \cq^{p-1-(r+k)} \left(\cp \left( \val{\cq^r \tilde f_k
  \sot \cq  ^{k-\ell+r} \tilde f_\ell  -\alpha ^{k- \ell+2r}
  f_{k,\ell,r} }  \right)\right)\\ 
&\leq  C\,  2^{r-\ell}\, \beta_r \alpha^{k-\ell+2r} \,
  \cq^{p-1-(r+k)}\left(\cp\left(g \otimes g\right)\right)\\ 
& \leq   C\, \beta_r \,  2^{-(k+\ell)/2}\, g_1^*,
\end{align*}
where we wrote (with $r'$ and $f$ replaced by $r$ and $f_k$ and by $k-\ell+r$ and $f_\ell$) that 
\begin{equation*}\label{eq:projvect}
\cq^{r'} \tilde{f} = \cq^{r'} \hat{f}  + \alpha^{r'} \sum_{j\in J} \theta_j^{r'}\crr_j (f)
\end{equation*}  
and used \reff{eq:unif-f-crit}, \reff{eq:cons-f-crit} and that
$(\beta_n,\, n\in \N)$ is non-decreasing for the second  inequality and
used \reff{eq:erg-bd} (with $f$ and $g$ replaced by $\cp\left(g \otimes g\right)$ and $g_1^*$) for the last. We deduce that:
\[
|H_{6,n} - \bar{H}_{6,n}| \leq \sum _{0\leq \ell<k\leq  p, \, r\geq 0}
|h_{k,\ell,r}^{(n)}-  \bar h_{k, \ell,r}^{(n)}| \,  \ind_{\{r+k< p\}}
\leq  C \left(\sum_{r=0}^n \beta_r \right)\, g_1^*. 
\] 
As $\lim_{n\rightarrow\infty } \beta_n=0$, we get that $\lim_{n\rightarrow \infty } n^{-1} \sum_{r=0}^n \beta_r =0$. We deduce from \reff{eq:lfgn-G} that a.s.:
\begin{equation}\label{eq:limH-tH}
\lim_{n\rightarrow \infty } n^{-1} |\G_{n-p}|^{-1} M_{\G_{n-p}} (|H_{6,n} -\bar{H}_{6,n}|) = 0.
\end{equation}
\medskip

We set $  H_6^{[n]} = \sum_{0\leq \ell<k\leq  p; \, r\geq 0}
  h_{k,\ell,r}\,  \ind_{\{r+k<  p\}} $ with  for $0\leq \ell<k\leq p$ and $0\leq r<p-k$: 
\begin{equation*}\label{eq:hat-crit}
 h_{k, \ell,r} = 2^{-(k+\ell)/2} \langle \mu, \cp f_{k, \ell,r} \rangle
=\langle \mu, \bar h_{k, \ell, r}^{(n)}  \rangle.
\end{equation*}
Notice that:
\begin{align*}
| \bar h_{k,\ell,r}^{(n)} -  h_{k,\ell,r}|
&\leq  2^{-(k+\ell)/2}
\sum_{j, j'\in J} \Big|\cq^{p-1-(r+k)} ( \cp (\crr_j f_k \sot \crr_{j'}
f_\ell)) - \langle \mu,\cp (\crr_j f_k \sot \crr_{j'}
f_\ell) \rangle \Big |\\
&\leq  2^{-(k+\ell)/2} \alpha^{p-1-(r+k)} \sum_{j, j'\in J}
  g_{k,\ell,j,j'}\\
&=  2^{-(k+\ell)/2} \alpha^{p-1-(r+k)} 
  g_{k,\ell},   
\end{align*}
where  we  used  \reff{eq:geom-erg}  (with   $f$  and  $g$  replaced  by
$\cp (\crr_j f_k \sot \crr_{j'} f_\ell)  $ and $g_{k,\ell,j, j'}$) for
the                 second                inequality                 and
$g_{k,\ell}=  \sum_{j, j'\in  J} g_{k,\ell,j,j'}$  for the  equality. We
have     that    $g_{k,     \ell}$     belongs     to    $F$.      Since
$|\cp f_{k,  \ell,r}|\leq \cp| f_{k,  \ell,r}|\leq 4 \cp  (g\otimes g)$,
thanks to the fourth inequality in \reff{eq:cons-f-crit}, we deduce from
\reff{eq:erg-bd} (with $f$ and $g$ replaced  by $4 \cp (g\otimes g)$ and
$g_2^*$) that for all $0\leq \ell<k$ and $ 0\leq r<p-k$:
\[
|\bar h_{k,\ell,r}^{(n)}|\leq  2^{-(k+\ell)/2} \,  g_2^* 
\quad\text{and}\quad
| h_{k,\ell,r}|\leq  2^{-(k+\ell)/2} \,   \langle \mu, g_2^*  \rangle.
\]
Set $r_0\in \N$ and $g_{r_0}=\sum_{0\leq
  \ell< k\leq  r_0} g_{k,\ell}$. Notice that $g_{r_0}$
belongs to $F$ and is non-negative. Furthermore, we have for $n$ large
enough so that $p> 2r_0$:
\begin{align*}
|\bar H_{6,n} - H_{6}^{[n]}| 
& \leq  \sum _{\substack{0\leq \ell<k\leq  p\\
 \, r\geq 0}}  | \bar h_{k,\ell,r}^{(n)}-  h_{k,\ell,r} | \,
  \ind_{\{r+k<  p\}} \\  
& \leq  \!\!\!\sum_{0\leq \ell< k\leq  r_0} \!\!\!\sum_{r=0}^{p-k-1}
  2^{-(k+\ell)/2} \alpha^{p-1 -(r+k)} g_{r_0} + \sum_{\substack{0\leq
  \ell< k\leq p\\   k>r_0}} \sum_{r=0}^{p-k-1} \left(|\bar
  h_{k,\ell,r}^{(n)}|+ | h_{k, \ell, r}|\right)    \,
  \ind_{\{r+k< p\}} \\ 
&\leq C  \,  g_{r_0}  + \sum_{0\leq  \ell< k\leq p,\,
  k>r_0} (p-k) \, 2^{ -(k+\ell)/2} \, (g_2^*+ \langle \mu, g_2^*
  \rangle ) \\ 
&\leq C \,  g_{r_0} + C n\,  2^{-r_0/2}  \,(g_2^*+ \langle \mu, g_2^* \rangle ). 
\end{align*}
We deduce that:
\[
\limsup_{n\rightarrow \infty } n^{-1} |\G_{n-p}|^{-1}\,
M_{\G_{n-p}}(|\bar H_{6,n} -  H_{6}^{[n]}|)\leq  C\;   2^{-r_0/2}
\langle \mu, g_2^* \rangle. 
\]
Since $r_0$ can be arbitrary large, we get that:
\begin{equation}\label{eq:limtH-hH2}
\lim_{n\rightarrow \infty } n^{-1} |\G_{n-p}|^{-1}\, M_{\G_{n-p}}(|\bar
H_{6,n} - H_{6}^{[n]}|) = 0. 
\end{equation}
\medskip

We set for $k, \ell\in \N$:
\begin{equation*}\label{eq:hklScri}
h_{k,\ell}^*= 2^{-(k+\ell)/2} \langle \mu, \cp (f^*_{k, \ell} ) \rangle.
\end{equation*}
Using the last inequality in \reff{eq:cons-f-crit} and the definition
\reff{eq:def-fkl*} of $f^*_{k, \ell}$, we deduce there exists a finite
constant $c$ independent of $n$ such that,  for all $k, \ell\in \N$, $|h_{ k,\ell}^*|\leq c  2^{-(k+\ell)/2}$. This
implies that $ H_0^*= \sum_{0\leq \ell <k} (k+1) |h_{k,\ell}^*| $ is finite and (see \reff{eq:S2-crit}) the sum: 
\begin{equation*}\label{eq:HstarV6}
H^{*}_{6}(\bF) =  \sum_{0\leq \ell <k} h_{k,\ell}^* = \scrit_2(\bF) 
\end{equation*}
is well defined and finite. We write:
\[
 h_{ k,\ell, r}= h_{k,\ell}^* +  h_{k,\ell, r}^{\circ},
\]
with 
\begin{equation*}\label{eq:hklnccri}
 h_{k,\ell,r}^{\circ} = 2^{-(k+\ell)/2} \langle \mu, \cp f^\circ_{k, \ell,r} \rangle,
\end{equation*} 
where we recall that $ f^\circ_{k, \ell, r}= f_{k, \ell, r}-
f^*_{k, \ell}$, 
and 
\begin{equation}\label{eq:decomphatHn}
H_{6}^{[n]}  = H_{6}^{[n],*} +H_{6}^{[n],\circ} 
\end{equation}
with  
\begin{equation*}\label{eq:Hnhat-H0nhat}
H_{6}^{[n],*} = \sum_{0\leq \ell<k\leq  p} (p-k) h^*_{k, \ell} \quad
\text{and} \quad 
H_{6}^{[n],\circ} = \sum_{0\leq \ell<k\leq  p; \, r\geq 0}  
h_{k,\ell,r}^\circ \,  \ind_{\{r+k<  p\}}.
\end{equation*} 
Recall $\lim_{n\rightarrow \infty } p/n=1$. We have:
\[
|n^{-1} H_{6}^{[n],*} - H^{*}_{6}(\bF)|\leq   |n^{-1} p - 1| |H^{*}_{6}(\bF)| + n^{-1} H_0^* +\sum_{\substack{0\leq  \ell<k \\ k>p}} |h_{k, \ell}^*| ,
\]
so that $\lim_{n\rightarrow \infty } |n^{-1}  H_{6}^{[n],*}- H^{*}_{6}(\bF)| = 0$ and
thus:
\begin{equation}\label{eq:limtH-hH}
\lim_{n\rightarrow \infty } n^{-1}  H_{6}^{[n], *} = H^{*}_{6}(\bF).
\end{equation}
We now prove that  $n^{-1}  H_{6}^{[n], \circ}$ converges towards 0. We have:
\begin{equation}\label{eq:detail-f0}
f^\circ_{k, \ell, r}= \sum_{j, j'\in J, \, \theta_j \theta_{j'}\neq 1} (\theta_{j'}\theta_{j})^r \theta_{j'}^{k-\ell} \,\, \crr_{j} f_k \sot \crr_{j'} f_\ell.
\end{equation}
This gives:
\begin{align}
| H^{[n],\circ}_{6}| 
& =\Big|  \sum_{0\leq \ell<k\leq  p, \, r\geq 0}  2^{ -(k+\ell)/2}
  \langle \mu, \cp f^\circ_{k, \ell, r}  \rangle\,  \ind_{\{r+k<
  p\}}\Big| \nonumber \\ 
& \leq  \sum_{0\leq \ell<k\leq  p}  2^{ -(k+\ell)/2}
  \sum_{j, j'\in J, \, \theta_j \theta_{j'}\neq 1} \Big| \langle \mu,
  \cp( \crr_{j} f_k \sot \crr_{j'} f_\ell) \rangle \Big| \,
  \Big|\sum_{r=0}^{p-k-1} (\theta_{j'}\theta_{j})^r
  \Big| \nonumber \\ 
& \leq c,   \nonumber
\end{align}
with $c=  c_J^2  \langle \mu, \cp (g\otimes g)  \rangle
  \sum_{0\leq \ell<k\leq  p}  2^{ -(k+\ell)/2} \sum_{j, j'\in J, \,
  \theta_j \theta_{j'}\neq 1}  |1- \theta_{j'}\theta_{j}|^{-1}$, 
and where we used \reff{eq:detail-f0} for the first inequality, the last
inequality of \reff{eq:cons-f-crit} for the second. Since  $J$ is
finite, we deduce that $c$ is finite. This gives that
$\lim_{n\rightarrow \infty } n^{-1}  H^{[n],\circ} _{6} = 0$. Recall that
$ H_{6}^{[n]}$ and $H^{*}_{6}(\bF)$ are complex numbers (\textit{i.e.} constant
functions). Use
\reff{eq:decomphatHn} and \reff{eq:limtH-hH} to get that: 
\begin{equation*}\label{eq:limhatHnV6}
\lim_{n\rightarrow \infty } n^{-1}  H_{6}^{[n]} = H^{*}_{6}(\bF) 
\end{equation*}
so that, as $ |\G_{n-p}|^{-1}\, M_{\G_{n-p}}(\cdot)$ is a probability measure,  a.s.:
\begin{equation}\label{eq:limtH-hH3}
\lim_{n\rightarrow \infty } n^{-1} |\G_{n-p}|^{-1}\, M_{\G_{n-p}}( H_{6}^{[n]})=H^{*}_{6}(\bF).
\end{equation}
\medskip

In conclusion, use \reff{eq:limH-tH}, \reff{eq:limtH-hH2}, \reff{eq:limtH-hH3}
and the definition \reff{eq:def-V6-crit} of $V_6(n)$ to deduce that a.s.:
\[
\lim_{n\rightarrow \infty } n^{-1}V_6(n)=H^{*}_{6}(\bF)= \sum_{0\leq
  \ell <k} 2^{-(k+\ell)/2} \langle \mu, \cp f_{k, \ell}^*   \rangle= \scrit_2(\bF) , 
\]
where $f_{k, \ell}^*$ is defined in \reff{eq:def-fkl*} and $\scrit_2(\bF) $
in \reff{eq:S2-crit}. Recall that:
\[
V_5(n)=  |\G_n|^{-1} \sum_{i\in \G_{n-p}} \sum_{0\leq \ell<k\leq  p } 2^{p-\ell} \cq^{p-k} \left( \tilde f_k \cq^{k-\ell} \tilde f_\ell\right)(X_i)= |\G_{n-p}|^{-1} M_{\G_{n-p}} (\Phi_n),
\]
where
\[
\Phi_n=\sum_{0\leq \ell<k\leq  p } 2^{-\ell} \cq^{p-k} \left( \tilde f_k
  \cq^{k-\ell} \tilde f_\ell\right).
\]
We have:
\[
|\Phi_n|
\leq  \sum_{0\leq \ell<k\leq  p } 2^{-\ell} \alpha ^{k-\ell} \cq^{p-k} (g^2)
\leq  \sum_{0\leq \ell<k\leq  p } 2^{-(k+\ell)/2} g_1\leq  C\,  g_1,
\]
where we used \reff{eq:geom-erg} for the first inequality and
\reff{eq:erg-bd} (with $f$ and $g$ replaced by $g^2$ and $g_1$) in the
second. Then, use \reff{eq:lfgn-G} to conclude that a.s.:
\[
\lim_{n\rightarrow\infty } n^{-1} V_5(n)=0.
\]
This ends the proof of the Lemma. 
\end{proof}

Using similar arguments as in the proof of Lemma \ref{lem:cvV2-crit}, we get
the following result. 
\begin{lem}\label{lem:cvV1-crit}
Under the assumptions of Theorem \ref{theo:critical}, we have that
a.s. $\lim_{n\rightarrow \infty } n^{-1}V_1(n) =\scrit_1(\bF)$ with
$\scrit_1(\bF)$ defined by \reff{eq:S1-crit} which is well defined and
finite.  
\end{lem}

\begin{proof}
We recall $V_1(n) = V_3(n) + V_4(n)$, see \reff{eq:DV4nsub} and thereafter for the
definition of $V_3(n)$ and $V_4(n)$. We first consider
the term $V_3(n)$. Recall that $V_3(n)=|\G_{n-p}|M_{\G{n-p}}(\Phi_n)$
with $\Phi_n=\sum_{\ell=0}^p 2^{-\ell}\,  \cq^{p-\ell} (\tilde
  f_\ell^2)$.
We have $\tilde f_\ell^2 \leq  g^2$ and $\cq^{p-\ell} (g^2)\leq g_1$ for
some $g_1\in F$ and 
thus $|\Phi_n|\leq  2 g_1$. We therefore deduce  that a.s. $\lim_{n\rightarrow \infty } n^{-1} V_3(n)=0$.
\medskip

We consider the term $V_4(n)=|\G_{n-p}|^{-1} M_{\G_{n-p}} (H_{4,n}) $
(see \reff{eq:V4}) with $H_{4,n}$ given by \reff{eq:def-H4n-sub-pt}:
\[
H_{4,n} = \sum_{\ell\geq 0, \, k\geq 0} h_{\ell,k}^{(n)}\,  \ind_{\{\ell+k<  p\}} \quad\text{and} \quad h_{\ell,k}^{(n)} = 2^{k-\ell}\,  \cq^{p-1-(\ell+k)} \left(\cp\left(\cq^k \tilde f_\ell \otimes^2 \right)\right).
\]
Recall $f_{\ell,\ell,k}$ defined in \reff{eq:def-fkl}. We set $\bar H_{4,n} = \sum_{\ell\geq 0, \, k\geq 0} \bar{h}_{\ell,k}^{(n)}\,  \ind_{\{\ell+k<  p\}}$ with
\[
\bar h_{\ell,  k}^{(n)} 
= 2^{k-\ell}\,\alpha^{2k}\,   \cq^{p-1-(\ell+k)} \left(\cp
  f_{\ell,\ell,k}\right)
= 2^{-\ell}\,   \cq^{p-1-(\ell+k)} \left(\cp f_{\ell,\ell,k}\right),
\]
where  we used  that  $2\alpha^2=1$. 
We have:
\begin{align*}
   |h_{\ell,k}^{(n)} - \bar h_{ \ell,k}^{(n)}|
&\leq  2^{k-\ell}
\cq^{p-1-(\ell+k)} \left(\cp \left( \val{\cq^k \tilde f_\ell \otimes
\cq  ^{k} \tilde f_\ell  -\alpha ^{2k}  f_{\ell,\ell,k } }  \right)\right)\\
&\leq  C\,  2^{k-\ell}\, \beta_k \, \alpha^{2k} \,  \cq^{p-1-(\ell+k)}\left(\cp\left(g \otimes
  g\right)\right)\\
& \leq  \beta_k\,  2^{-\ell}\, g_1^*,
\end{align*}
with $g_1^*\in F$, where we used 
 \reff{eq:unif-f-crit},
with the representation $\cq^k \tilde{f}_\ell = \cq^k \hat{f}_\ell  +
\alpha^k \sum_{j\in J} \theta_j^k \crr_j (f_\ell)$, 
  \reff{eq:cons-f-crit}  for the
second  inequality and \reff{eq:erg-bd}  for the last.
We deduce  that:
\[
|H_{4,n} - \bar  H_{4,n}| \leq \sum _{\ell\geq 0, \, k\geq 0}
|h_{\ell,k}^{(n)}-  \bar h_{ \ell,k}^{(n)}| \,  \ind_{\{\ell+k<p\}} \leq  2 \left(\sum_{k=0}^n \beta_k \right)\, g_1^*.
\] 
As $\lim_{n\rightarrow\infty } \beta_n=0$, we get that
$\lim_{n\rightarrow \infty } n^{-1} \sum_{k=0}^n \beta_k =0$. 
We deduce from \reff{eq:lfgn-G} that a.s.:
\begin{equation}\label{eq:limH-tH-4}
\lim_{n\rightarrow \infty } n^{-1} |\G_{n-p}|^{-1} M_{\G_{n-p}} (|H_{4,n} - \bar H_{4,n}|)=0.
\end{equation}
\medskip

We set $ H_{4}^{[n]} = \sum_{\ell\geq 0, \, k\geq 0}  
h_{\ell,k}\,  \ind_{\{\ell+k<  p\}}$ with:
\[
 h_{\ell,k} = 2^{-\ell} \langle \mu, \cp f_{\ell, \ell,k} \rangle.
\]
Notice that:
\begin{align*}
| \bar h_{\ell,k}^{(n)} -  h_{\ell,k}|
&\leq  2^{-\ell}
\sum_{j, j'\in J} \Big|\cq^{p-1-(\ell+k)} ( \cp (\crr_j f_\ell\sot \crr_{j'}
f_\ell)) - \langle \mu,\cp (\crr_j f_\ell\sot \crr_{j'}
f_\ell) \rangle \Big |\\
&\leq  2^{-\ell} \alpha^{p-1-(\ell+k)} \sum_{j, j'\in J}
  g_{\ell,j,j'}\\
&=  2^{-\ell} \alpha^{p-1-(\ell+k)}
  g_{\ell},   
\end{align*}
where  we  used  \reff{eq:geom-erg}  (with   $f$  and  $g$  replaced  by
$\cp (\crr_j f_\ell\otimes \crr_{j'} f_\ell)  $ and $g_{\ell,j, j'}$) for
the                 second                inequality                 and
$g_{\ell}=  \sum_{j, j'\in  J} g_{\ell,j,j'}$  for the  equality. We
have     that    $g_{     \ell}$     belongs     to    $F$.      Since
$|\cp f_{\ell,  \ell,k}|\leq \cp| f_{\ell,  \ell,k}|\leq 4 \cp  (g\otimes g)$,
thanks to the fourth inequality in \reff{eq:cons-f-crit}, we deduce from
\reff{eq:erg-bd} (with $f$ and $g$ replaced  by $4 \cp (g\otimes g)$ and
$g_2^*$) that:
\[
|\bar h_{\ell,k}^{(n)}|\leq  2^{-\ell} \,  g_2^* 
\quad\text{and}\quad
|h_{\ell,k}|\leq  2^{-\ell} \,   \langle \mu, g_2^*  \rangle.
\]
Set $r_0\in \N$ and $g_{r_0}=\sum_{0\leq
  \ell\leq  r_0} g_{\ell}$. Notice that $g_{r_0}$
belongs to $F$ and is non-negative. Furthermore, we have for $n$ large
enough so that $p> 2r_0$:
\begin{align*}
|\bar H_{4,n} - H_{4}^{[n]}| 
&\leq  \sum _{\ell\geq 0, \,  k\geq 0 }  | \bar  h_{\ell,k}^{(n)}- 
  h_{\ell,k} | \,  \ind_{\{\ell+k<  p\}} \\ 
&\leq  \sum_{0\leq \ell\leq  r_0, \, k\geq 0}  2^{-\ell} \alpha^{p-1
  -(\ell+k)} g_{r_0}\ind_{\{\ell +k <p\}}
 + \sum_{\ell>r_0, \, k\geq 0} \left(|\bar
  h_{\ell,k}^{(n)}|+ | h_{ \ell,k}|\right)    \,
  \ind_{\{\ell+k<    p\}}\\ 
&\leq C  \,  g_{r_0} + \sum_{\ell>r_0} (p-\ell) \, 2^{-\ell} \, (g_2^*+
  \langle \mu, g_2^* \rangle )\ind_{\{\ell<p\}} \\ 
&\leq C \,  g_{r_0} + n\,  2^{-r_0} \,(g_2^*+\langle \mu, g_2^* \rangle ). 
\end{align*}
We deduce that:
\[
\limsup_{n\rightarrow \infty } n^{-1} |\G_{n-p}|^{-1}\,
M_{\G_{n-p}}(|\bar  H_{4,n} -  H_{4}^{[n]}|)\leq    2^{
1-r_0} \langle \mu, g_2^* \rangle.
\]
Since $r_0$ can be arbitrary large, we get that a.s.:
\begin{equation}\label{eq:limtH-hH2-4}
\lim_{n\rightarrow \infty } n^{-1} |\G_{n-p}|^{-1}\,
M_{\G_{n-p}}(|\bar  H_{4,n} -  H_{4}^{[n]}|) = 0. 
\end{equation}
\medskip

Now, we study the limit of $H_4^{[n]}$. 
We set for $k, \ell\in \N$:
\[
h_{\ell}^*= 2^{-\ell} \langle \mu, \cp f^*_{\ell, \ell} \rangle.
\]
Using the last inequality in \reff{eq:cons-f-crit} and the definition
\reff{eq:def-fkl*} of $f^*_{\ell, \ell}$, we deduce there exists a finite
constant $c$ independent of $n$ (but depending on $\bF$) such that,  for all $ \ell\in \N$,
$|h_{\ell}^*|\leq c  2^{-\ell}$. This 
implies that $H_0^*= \sum_{\ell\geq 0 } (\ell+1) |h_{\ell}^*|$ is finite
and the sum 
\begin{equation*}\label{eq:HstarV4}
H^{*}_{4}(\bF)= \sum_{ \ell\geq 0} h_{\ell}^* 
\end{equation*}
is well defined and finite. 
We write:
\[
 h_{ \ell,k}= h_{\ell}^* +  h_{\ell,k}^{\circ},
\]
with $ h_{\ell,k}^{\circ}=
 2^{-\ell} \langle \mu, \cp f^\circ_{\ell, \ell,k} \rangle$, where $
 f^\circ_{\ell, \ell,k}=f_{\ell, \ell,k} - f^*_{\ell, \ell}$ is defined
 in \reff{eq:def-fkl},   
and 
\begin{equation*}\label{eq:decomphatHn-4}
H_4^{[n]} =  H^{[n],*}_{4}  +  H^{[n],\circ}_{4},
\end{equation*}
with  $  H_{4}^{[n], *}=\sum_{\ell\geq 0} (p-\ell) h^*_{ \ell}$ and
$ H_{4}^{[n], \circ} = \sum_{\ell\geq 0, k\geq 0}  
h_{\ell,k}^{\circ}\,  \ind_{\{\ell+k<  p\}}$. 
We have:
\[
|n^{-1}  H^{[n],*}_{4} - H^{*}_{4}(\bF)| \leq   |n^{-1} p - 1| H^{*}_{4}(\bF) + n^{-1} H_0^* +\sum_{\ell >p} |h_{ \ell}^*| ,
\]
so that $\lim_{n\rightarrow \infty } |n^{-1}   H^{[n],*}_{4} - H^{*}_{4}(\bF)| = 0$ and thus:
\begin{equation}\label{eq:limtH-hH-4}
\lim_{n\rightarrow \infty } n^{-1}   H^{[n],*}_{4}=H^{*}_{4}(\bF).
\end{equation}
We now prove that  $n^{-1} H_{4}^{[n], \circ}$ converges towards 0. 
We have:
\begin{equation}\label{eq:detail-f0-4}
f^\circ_{\ell, \ell, k}= \sum_{j, j'\in J, \, \theta_j \theta_{j'}\neq 1} (\theta_{j'}\theta_{j})^k  \,\, \crr_{j'} f_\ell \sot \crr_j f_\ell.
\end{equation}
This gives:
\begin{align*}
| H^{[n], \circ}_{4}| 
&=\Big|  \sum_{\ell\geq 0, k\geq 0}  2^{-\ell} \langle \mu, \cp
  f^\circ_{\ell, \ell, k}  \rangle\,  \ind_{\{\ell+k< p\}}\Big|\\ 
&\leq  \sum_{\ell\geq 0}  2^{-\ell} \sum_{j, j'\in J, \,
  \theta_j \theta_{j'}\neq 1}  \Big| \langle \mu, \cp(
  \crr_{j'} f_\ell \otimes \crr_j f_\ell) \rangle \Big| \,
  \Big|\sum_{k=0}^{p-\ell-1} (\theta_{j'}\theta_{j})^k \Big|\\ 
&\leq c,
\end{align*}
with $c= c_J^2  \langle \mu, \cp (g\otimes g)  \rangle \sum_{\ell\geq 0}
  2^{-\ell}\sum_{j, j'\in J, \, \theta_j \theta_{j'}\neq 1}
{|1- \theta_{j'}\theta_{j}|}^{-1}$, and 
where we used  \reff{eq:detail-f0-4} for the first  inequality, the last
inequality of \reff{eq:cons-f-crit} for the second. Since $J$ is finite,
we    deduce     that  $c$     is    finite.    This     gives    that
$\lim_{n\rightarrow \infty  } n^{-1}   H^{[n],\circ}  _{4} =  0$. Recall
that   $  H_{4}^{[n]}$   and   $H^{*}_{4}(\bF)$   are  complex   numbers
(\textit{i.e.}  constant  functions).  Use \reff{eq:limtH-hH-4}  to  get
that:
\begin{equation*}\label{eq:limhatHnV4}
\lim_{n\rightarrow \infty } n^{-1}  H_{4}^{[n]}=H^{*}_{4}(\bF) 
\end{equation*}
so that, as $ |\G_{n-p}|^{-1}\, M_{\G_{n-p}}(\cdot)$ is a probability measure,  a.s.:
\begin{equation}\label{eq:limtH-hH3-4}
\lim_{n\rightarrow \infty } n^{-1} |\G_{n-p}|^{-1}\, M_{\G_{n-p}}(
H_{4}^{[n]})=H^{*}_{4}(\bF). 
\end{equation}
\medskip

In conclusion, use \reff{eq:limH-tH-4}, \reff{eq:limtH-hH2-4}, \reff{eq:limtH-hH3-4}
and the definition \reff{eq:V4} of $V_4(n)$ to deduce that a.s.:
\[
\lim_{n\rightarrow \infty } n^{-1}V_4(n)=H^{*}_{4}(\bF)= \sum_{\ell\geq 0}
 2^{-\ell} \langle \mu, \cp f_{\ell, \ell}^*  \rangle=\scrit_1(\bF) ,
\]
where $f_{\ell, \ell}^*$ is defined in \reff{eq:def-fkl*} and
$\scrit_1(\bF) $ in \reff{eq:S1-crit}. 
\end{proof}

 The proof of the next Lemma is a direct consequence of \reff{eq:def-V} and
 Lemmas \ref{lem:cvR2-crit}, \ref{lem:cvV1-crit} and \ref{lem:cvV2-crit}.
\begin{lem}
   \label{lem:cvV-crit}
Under the assumptions of Theorem \ref{theo:critical}, we have the
following convergence in probability:
\[
\lim_{n\rightarrow\infty } n^{-1} V(n)=\scrit_1(\bF)+ 2\scrit_2(\bF), 
\]
where $\scrit_1(\bF)$ and  $\scrit_2(\bF)$, defined by \reff{eq:S1-crit}
and \reff{eq:S2-crit}, are well defined and finite.
\end{lem}

We now check the Lindeberg condition. Recall $R_3(n)$ defined  in \reff{eq:def-R3}. 

\begin{lem}
   \label{lem:cvG-crit}
Under the assumptions of Theorem \ref{theo:critical}, we have that
$\lim_{n\rightarrow\infty } n^{-2} R_3(n)=0$.
\end{lem}
\begin{proof}
Keeping the notation of Lemma \ref{lem:cvG}, using Lemma \ref{lem:M4}
(with the main contribution coming from $\psi_{8,n}$ and $\psi_{9,n}$
therein), we get (compare with \reff{eq:majo-hnl-g1}) that 
 for 
$n\geq p \geq \ell\geq 0$:
\[
|h_{n,\ell}|\leq  (p-\ell)^2 2^{2(p-\ell)} g_1,
\]
with $h_{n,\ell}(x)=\E_x\left[M_{\G_{p-\ell}}(\tilde f_\ell) ^4\right]$
and  $g_1\in F$. 
Following the proof of Lemma \ref{lem:cvG}, we get that:
\[
n^{-2} R_3(n) \leq    16 \, n^3\, 2^{-(n-p)} \, \langle \nu, g_2 \rangle.
\]
This  ends the proof as $\lim_{n\rightarrow \infty } p=\infty $
and $\lim_{n\rightarrow\infty } n-p - \lambda \log(p)=+\infty $ for all $\lambda>0$.
\end{proof}

The proof of  Theorem \ref{theo:critical} mimics  then the proof of
Theorem
\ref{theo:subcritical}.

\section{Proof of Lemma \ref{lem:martingale} and of Theorem
  \ref{theo:super-critical}} 
\label{sec:proof-super}

\subsection{Proof of Lemma \ref{lem:martingale}}
\label{sec:martingaleSuC}

Let $f\in F$ and $j\in J$. Use that $\crr_j(F)\subset \C F$ to deduce that
$\E\left[|M_{n,j}(f)|^2\right]$ is finite. 
We have for $n\in \N^*$:
\begin{align*}
\E[M_{n,j}(f)|\Hh_{n-1}] 
&= (2\alpha_j)^{-n} \sum_{i\in\GG_{n-1}} \E[\crr_j f(X_{i0}) +
  \crr_j f(X_{i1})|\Hh_{n-1}]\\ 
&= (2\alpha_j)^{-n} \sum_{i\in\GG_{n-1}} 2\, \cq \crr_j f(X_{i}) \\ 
&= (2\alpha_j)^{-(n-1)}\sum_{i\in\GG_{n-1}} \crr_j f (X_{i}) \\
& = M_{n-1,j}(f),
\end{align*}
where the second equality follows from branching Markov property and the
third  follows from the fact that $\crr_j$ is the projection on the
eigen-space associated to the eigen-value $\alpha_j$ of $\cq$. This
gives that $M_{j}(f)$ is a $\ch$-martingale. 
We also  have, writing $f_j$ for $\crr_j(f)$:
\begin{align}
\nonumber
\E\left[|M_{n,j}(f)|^2\right]
&= (2\alpha)^{-2n} \, \E\left[M_{\G_n} (f_j) M_{\G_n} (\overline
  f_j)\right]\\
\nonumber
&= (2\alpha^2)^{-n} \,\langle \nu, \cq^n (|f_j|^2)  \rangle
+ (2\alpha)^{-2n} \sum_{k=0}^{n-1} 2^{n+k}\langle \nu, \cq^{n-k-1}
  \cp\left(\cq^k f_j  \sot \cq^{k} \overline f_j\right) \rangle\\
\label{eq:mart-majo-L2}
&= (2\alpha^2)^{-n} \,\langle \nu, \cq^n (|f_j|^2)  \rangle
+ (2\alpha^2)^{-n} \sum_{k=0}^{n-1} (2\alpha^2)^{k} \langle \nu, \cq^{n-k-1}
  \cp\left( f_j  \sot  \overline f_j\right) \rangle\\
\nonumber
&\leq (2\alpha^{2})^{-n}\langle\nu,g_{1}\rangle +
  (2\alpha^{2})^{-n}\sum_{k=0}^{n-1}
(2\alpha^{2})^{k}  \langle\nu,g_{2}\rangle \\ 
\nonumber
&\leq \langle\nu,g_{3}\rangle,
\end{align}
where  we used  the  definition  of $M_{n,j}$  for  the first  equality,
\eqref{eq:Q2-bis} with $m=n$ for the  second equality, the fact that  $f_j$
(resp. $\overline f_j$)  is an
eigen-function  associated  to  the  eigenvalue $\alpha_j$
(resp. $\overline \alpha_j$) for  the  third
equality, \eqref{eq:erg-bd} twice (with $f$ and $g$ replaced by
$|f_j|^2$ and $g_{1}$ and by $\cp\left( f_j  \sot  \overline
  f_j\right) $ and $g_2$) for the first inequality  and  $2\alpha^{2}  >  1$ as well as $g_3=g_1+
g_2/(2\alpha^2 - 1)$ for  the  last
inequality. Since $g_3$ belongs to $F$ and does   not  depend on  $n$,  this
implies     that
$\sup_{n\in\N}\E\left[|M_{n,j}(f)|^2\right] < + \infty.$ Thus the
martingale  $M_{j}(f)$ converges a.s. and  in  $L^{2}$  towards a limit.

\subsection{Proof of Theorem \ref{theo:super-critical}}
\label{sec:proof-super-thm}
Recall the sequence $(\beta_{n},n\in\NN)$ defined in Assumption
\ref{hyp:F3} and  the $\sigma$-field $\Hh_{n} = \sigma\{X_{u},u\in\TT_{n}\}$.
Let  $(\p_{n},n\in\NN)$ be a sequence of integers such that $\p_n$ is even
and (for $n\geq 3$):
\begin{equation}\label{eq:hypSCrit-p}
\frac{5n}{6} <\p_n<n, \quad  \lim_{n\rightarrow \infty } (n-\p_n)=\infty
\quad  \text{and} \quad \lim_{n\rightarrow \infty} 
\alpha^{-(n-\p_n)} \beta_{\p_{n}/2} = 0.  
\end{equation}
Notice such sequences exist. 
When there is no ambiguity, we shall write $\p$  for $\p_n$.
 We deduce from \eqref{eq:nof-D} that:
\begin{equation}\label{eq:decomSCrit}
N_{n, \emptyset}(\bF) = R_{0}(n) + R_{4}(n) + T_{n}(\bF),
\end{equation}
with notations from \reff{eq:N=D+R} and   \reff{eq:reste01}:
\begin{align*}
R_{0}(n) &= |\GG_{n}|^{-1/2} \sum_{k=0}^{n-\p_n-1}
  M_{\GG_{k}}(\tilde{f}_{n-k}),\\
 T_{n}(\bF)=R_1(n) &= \sum_{i\in\GG_{n-\p_n}} \EE[N_{n,i}(\bF)|\Hh_{n-\p_n}],\\
 R_{4}(n)=\Delta_n &= \sum_{i\in \GG_{n-\p_n}} \left(N_{n,i}(\bF) -
            \EE[N_{n,i}(\bF)|\Hh_{n-\p_n}]\right). 
\end{align*}
Furthermore, using the branching Markov property, we get for all $i\in \G_{n-\p_n}$:
\begin{equation}\label{eq:bmpSCrit}
\EE[N_{n,i}(\bF)|\Hh_{n-\p_n}] = \EE[N_{n,i}(\bF)|X_{i}].
\end{equation}
We have the following elementary lemma.
\begin{lem}\label{lem:cvR0SCrit}
Under the assumptions of Theorem \ref{theo:super-critical}, we have the
following convergence: 
\begin{equation*}
\lim_{n\rightarrow \infty} (2\alpha^{2})^{-n}\, \EE\left[R_{0}(n)^{2}\right] = 0.
\end{equation*}
\end{lem}
\begin{proof}
We follow the proof of Lemma \ref{lem:cvR0}. As $2\alpha^{2} > 1$, we
get that
$\EE[M_{\GG_{k}}(\tilde{f}_{n-k})^{2}] \leq
2^{k}(2\alpha^{2})^{k}\langle \nu,g\rangle$ for  some $g\in F$ and all
$n\geq k\geq 0$. This 
implies, see \reff{eq:majo-R0}, that for some constant $C$ which does
not depend on $n$ or $\p$:
\begin{equation*}
\EE\left[R_{0}(n)^{2}\right]^{1/2} \leq C\, 
2^{-\p/2} (2\alpha^2)^{(n-\p)/2} .
\end{equation*}
It follows from the previous inequality that $
(2\alpha^{2})^{-n} \E\left[R_{0}(n)^{2}\right] \leq 
C  (2\alpha)^{-2\p}$. 
Then use $2\alpha > 1$ and $\lim_{n\rightarrow \infty} \p = \infty$ to conclude.
\end{proof}
Next, we have the following lemma.
\begin{lem}\label{lem:cvR4SCrit}
Under the assumptions of Theorem \ref{theo:super-critical}, we have the following convergence:
\begin{equation*}
\lim_{n\rightarrow \infty} (2\alpha^{2})^{-n}\EE\left[R_{4}(n)^{2}\right] = 0.
\end{equation*}
\end{lem}
\begin{proof}
First, we have:
\begin{align}
\EE[R_{4}(n)^{2}] 
&= \EE\left[\left(\sum_{i\in \GG_{n-\p}} (N_{n,i}(\bF) -
  \EE[N_{n,i}(\bF)|X_{i}])\right)^{2}\right] \nonumber\\
 &= \EE\left[\sum_{i\in \GG_{n-\p}}\EE[(N_{n,i}(\bF) -
   \EE[N_{n,i}(\bF)|X_{i}])^{2}|\Hh_{n-\p}]\right]
 \nonumber \\ 
&\leq  \EE\left[\sum_{i\in \GG_{n-\p}}\EE[N_{n,i}(\bF)^2|X_i]\right],
\label{eq:R4SCrit1} 
\end{align}
where we used \eqref{eq:bmpSCrit} for the first equality and the branching Markov chain property for the second and the last inequality.  Note that for all $i\in\GG_{n-\p}$ we have 
\begin{align*}
\EE\left[\EE[N_{n,i}(\bF)^2|X_i]\right]
&= |\GG_{n}|^{-1}  \EE\left[\left(\sum_{\ell=0}^{\p}
  M_{i\GG_{\p-k}}(\tilde{f}_{\ell})\right)^{2}|X_{i}\right]
\\
&\leq |\GG_{n}|^{-1} \left(\sum_{\ell=0}^{\p}
  \EE_{X_{i}}[M_{\GG_{\p-\ell}}(\tilde{f}_{\ell})^{2}]^{1/2}\right)^{2},
\end{align*}
where we used the definition of $N_{n,i}(\bF)$ for the equality and the Minkowski's inequality for the last inequality. We have: 
\begin{align*}
\EE_{X_{i}}[M_{\GG_{\p-\ell}}(\tilde{f}_{\ell})^{2}] 
&= 2^{\p-\ell} \Qq^{\p-\ell}(\tilde{f}_{\ell}^{2})(X_{i}) + \sum_{k=0}^{\p-\ell-1}
  2^{\p-\ell+k} \Qq^{\p-\ell-k -1}(\Pp(\Qq^{k}\tilde{f}_{\ell} \otimes
  \Qq^{k}\tilde{f}_{\ell}))(X_{i})
 \\
&\leq 2^{\p-\ell} g_{2}(X_{i}) + \sum_{k=0}^{\p-\ell-1}
  2^{\p-\ell+k} \alpha^{2k}\Qq^{\p-\ell-k-1}(\Pp(g_{1}\otimes
  g_{1}))(X_{i})
 \\
&\leq 2^{\p-\ell} g_{2}(X_{i}) + \sum_{k=0}^{\p-\ell-1}
  2^{\p-\ell}(2\alpha^{2})^{k}g_{3}(X_{i}) \\
&\leq (2\alpha)^{2(\p-\ell)}g_{4}(X_{i}),
\end{align*}
where we used \eqref{eq:Q2} for the first equality, (ii) of Assumption
\ref{hyp:F} and \eqref{eq:unif-f} for the first inequality,
\eqref{eq:erg-bd} and (iv) of Assumption \ref{hyp:F} for the second, and
 $2\alpha^{2} > 1$  for the last. The latter inequality implies that,
 with $g_5$ equal to $g_4$ up to a finite multiplicative constant:
\begin{equation}
 \label{eq:R4SCrit4}
\EE[N_{n,i}(\bF)^{2}|X_{i}]] 
\leq |\GG_{n}|^{-1}\left(\sum_{\ell=0}^{\p}
  (2\alpha)^{(\p-\ell)}\right)^{2}g_{4}(X_{i})
  =2^{-n} (2\alpha)^{2\p} \,  g_{5}(X_{i}).
\end{equation}
Using \eqref{eq:R4SCrit1}, \eqref{eq:R4SCrit4} and   \eqref{eq:Q1} as
well as Assumption \ref{hyp:F1} with $g_6\in F$, we get:
\[
(2\alpha^{2})^{-n}\EE[R_{4}(n)^{2}] 
\leq (2\alpha^{2})^{-n} 2^{n-\p} 2^{-n} (2\alpha)^{2\p} \langle \nu,
  \cq^n g_5 \rangle
\leq  (2\alpha^2)^{-(n-\p)} \langle \nu,
   g_6 \rangle.
\]
We then conclude using  that $2\alpha^{2} > 1$ and $\lim_{n\rightarrow
  \infty }(n- \p)=\infty $.
\end{proof}

Now, we study the third term of the right hand side of
\eqref{eq:decomSCrit}. First, note that: 
\begin{align*}
T_{n}(\bF) 
&= \sum_{i\in\GG_{n-\p}}  \EE[N_{n,i}(\bF)|X_{i}] \\ 
&= \sum_{i\in\GG_{n-\p}}|\GG_{n}|^{-1/2} \sum_{\ell = 0}^{\p}
  \EE_{X_{i}}[M_{\GG_{\p - \ell}}(\tilde{f}_{\ell})] \\ 
&= |\GG_{n}|^{-1/2} \sum_{i\in\GG_{n-\p}} \sum_{\ell = 0}^{\p} 2^{\p -
  \ell} \Qq^{\p - \ell}(\tilde{f}_{\ell})(X_{i}), 
\end{align*}
where we used \eqref{eq:bmpSCrit} for the first equality, the definition
\reff{eq:def-NiF}  of  $N_{n}(\bF)$  for  the  second  equality  and
\eqref{eq:Q1} for the last equality. Next, projecting in the eigen-space
associated to the eigenvalue $\alpha_{j}$, we get
\begin{equation*}\label{eq:decomTSCrit}
T_{n}(\bF) = T^{(1)}_{n}(\bF) + T^{(2)}_{n}(\bF),
\end{equation*}
where, with $\hat f= f - \langle \mu, f \rangle - \sum_{j\in J}
\crr_j(f)$ defined in \reff{eq:fhatcrit-S}:
\begin{align*}
T^{(1)}_{n}(\bF) 
&= |\GG_{n}|^{-1/2} \sum_{i\in \GG_{n-\p}} \sum_{\ell=0}^{\p} 2^{\p-\ell}
  \left(\Qq^{\p-\ell}(\hat{f}_{\ell}) \right)(X_{i}),\\ 
T^{(2)}_{n}(\bF) 
& = |\GG_{n}|^{-1/2} \sum_{i\in \GG_{n-\p}}
  \sum_{\ell=0}^{\p} 2^{\p-\ell} 
 \alpha^{\p-\ell}\sum _{j\in J} 
\theta_j^{\p-\ell}\crr_j (f_\ell)
(X_{i}). 
\end{align*}
We have the following lemma.
\begin{lem}\label{lem:cvT1SCrit}
Under the assumptions of Theorem \ref{theo:super-critical}, we have the following convergence:
\begin{equation*}
\lim_{n\rightarrow \infty}(2\alpha^{2})^{-n/2} \EE[|T^{(1)}_{n}(\bF)|] = 0.
\end{equation*}
\end{lem}
\begin{proof}
We have:
\begin{align*}
(2\alpha^{2})^{-n/2}\EE[|T^{(1)}_{n}(\bF)|] 
&\leq (2\alpha)^{-n} \EE\left[\sum_{i\in \GG_{n-\p}} \sum_{\ell=0}^{\p} 2^{\p-\ell}
  |\Qq^{\p-\ell}(\hat {f}_{\ell})(X_i)  |\right] \\ 
&\leq (2\alpha)^{-n} \EE\left[\sum_{i\in\GG_{n-\p}}\sum_{\ell=0}^{\p} 2^{\p-\ell}
 \alpha^{\p-\ell}  \beta_{\p-\ell}g(X_{i})\right] \\ 
&=\sum_{\ell=0}^\p 2^{-\ell}\alpha^{-(n-\p+\ell)}  \beta_{\p-\ell} \,
  \langle  \nu, \cq^{n-\p} g \rangle,
\end{align*}
where we used the definition of $T^{(1)}_{n}(\bF)$ for the first
inequality, \eqref{eq:unif-f-crit} for the second  and \eqref{eq:Q1} for
the  equality. 
Using \eqref{eq:erg-bd} and the property (iii), we get that there is a
finite  positive constant $C$ independent of $n$ and $\p$ such that
$\langle\nu,\Qq^{n-\p}g\rangle < C$. 
We have:
\[
\sum_{\ell=0}^{\p/2} 2^{-\ell}\alpha^{-(n-\p+\ell)}  \beta_{\p-\ell} 
\leq \alpha^{-(n-\p)}  \beta_{\p/2}  \sum_{\ell=0}^{\p/2} (2\alpha)^{-\ell}. 
\]
Using the third condition in \reff{eq:hypSCrit-p} and that $2\alpha>1$, we deduce the right hand-side converges to $0$ as $n$ goes to infinity. Without loss of generality, we can assume that the sequence $(\beta_n, n\in \N^*)$ is bounded by 1. Since $\alpha>1/\sqrt{2}$, we also have:
\[
\sum_{\ell=\p/2}^\p 2^{-\ell}\alpha^{-(n-\p+\ell)}  \beta_{\p-\ell} \leq  (1- 2\alpha)^{-1}  \, 2^{-\p/2} \alpha^{ -n + \p/2} \leq (1- 2\alpha)^{-1}  \, 2^{ n/2 - 3\p/4}.
\]
Using that $n/2 - 3\p/4 < -n/8$, thanks to  the first condition in \reff{eq:hypSCrit-p},  we deduce the right hand-side converges to $0$ as $n$ goes to infinity.
Thus,   we get that  $ \lim_{n \rightarrow \infty}(2\alpha^{2})^{-n/2} \EE[|T^{(1)}_{n}(\bF)|] = 0$.
\end{proof}

Now, we deal with the term $T^{(2)}_{n}(\bF)$ in the following
result. Recall $M_{\infty , j}$ defined in Lemma \ref{lem:martingale}. 
\begin{lem}\label{lem:cvT2SCrit}
Under the assumptions of Theorem \ref{theo:super-critical}, we have the
following convergence: 
\begin{equation*}
(2\alpha^{2})^{-n/2}T^{(2)}_{n}(\bF) - \sum_{\ell\in \N} (2\alpha)^{-\ell} \sum_{j\in J} \theta_j^{n-\ell}M_{\infty  ,j}(f_\ell)  \; \xrightarrow[n\rightarrow \infty ]{\P}  \; 0.
\end{equation*}
\end{lem}
\begin{proof}
By definition of $T^{2}_{n}(\bF)$, we have
$T^{2}_{n}(\bF)=2^{-n/2} \sum_{\ell = 0}^{\p} (2\alpha)^{n-\ell} \sum_{j\in J}
\theta_{j}^{n-\ell} M_{n,j}(f_{\ell})$ and thus:
\begin{multline}\label{eq:T2-0}
(2\alpha^{2})^{-n/2}T^{(2)}_{n}(\bF) - \sum_{\ell\in \N}
(2\alpha)^{-\ell} \sum_{j\in J} \theta_j^{n-\ell}M_{\infty  ,j}(f_\ell)
\\ 
= \sum_{\ell = 0}^{\p} (2\alpha)^{-\ell} \sum_{j\in J}
\theta_{j}^{n-\ell} (M_{n,j}(f_{\ell}) - M_{\infty,j}(f_{\ell})) -
\sum_{\ell = \p+1}^{\infty} (2\alpha)^{-\ell} \sum_{j\in J}
\theta_{j}^{n-\ell} M_{\infty,j}(f_{\ell}). 
\end{multline}
Using that $|\theta_j|=1$, we get:
\[
\EE[|\sum_{\ell = 0}^{\p} (2\alpha)^{-\ell} \sum_{j\in J}
\theta_{j}^{n+\p-\ell} (M_{n,j}(f_{\ell}) - M_{\infty,j}(f_{\ell}))| ] 
\leq \sum_{\ell = 0}^{\p} (2\alpha)^{-\ell} \sum_{j\in J}
\EE[|M_{n,j}(f_{\ell}) - M_{\infty,j}(f_{\ell})|].
\]
 Now, using \eqref{eq:unif-f}, a close inspection of the proof of Lemma
 \ref{lem:martingale}, see \reff{eq:mart-majo-L2}, 
 reveals us that there exists a finite constant $C$ (depending on $\bF$)
 such that  for all $j\in J$, we have:
\begin{equation*}
\label{eq:supMn-bd}
\sup_{\ell\in\NN}\sup_{n\in\NN}\EE[|M_{n,j}(f_{\ell})|^{2}] \leq C.
\end{equation*}
The $L^2(\nu)$ convergence in Lemma \ref{lem:martingale} yields that:
\begin{equation}
\label{eq:supMi-bd}
\sup_{\ell\in\NN}\EE[|M_{\infty ,j}(f_{\ell})|^{2}] \leq C.
\quad\text{and}\quad
\sup_{\ell \in \NN} \sup_{n \in \NN} \sum_{j\in J}\EE[|M_{n,j}(f_{\ell}) - M_{\infty,j}(f_{\ell})|] < 2|J|\sqrt{C}.
\end{equation}
Since Lemma \ref{lem:martingale} implies that $\lim_{n\rightarrow \infty
} \EE[|M_{n,j}(f_{\ell}) - M_{\infty,j}(f_{\ell})|]=0$, we deduce, as $2\alpha>1$ by the
dominated convergence theorem that:
\begin{equation}\label{eq:T2-1}
\lim_{n\rightarrow +\infty}\EE[|\sum_{\ell = 0}^{\p} (2\alpha)^{-\ell} \sum_{j\in J} \theta_{j}^{n+\p-\ell} (M_{n,j}(f_{\ell}) - M_{\infty,j}(f_{\ell}))|] = 0.
\end{equation}

On the other hand, we have
\begin{equation}\label{eq:T2-2}
\EE[|\sum_{\ell = \p+1}^{\infty} (2\alpha)^{-\ell} \sum_{j\in J}
  \theta_{j}^{n-\ell} M_{\infty,j}(f_{\ell})|] 
\leq \sum_{\ell = \p+1}^{\infty} (2\alpha)^{-\ell} \sum_{j\in J}
  \EE[|M_{\infty,j}(f_{\ell})|] 
 \leq |J|\sqrt{C}\, \sum_{\ell = \p+1}^{\infty} (2\alpha)^{-\ell} , 
\end{equation}
where we used $|\theta_{j}| = 1$ for the first inequality and the
Cauchy-Schwarz inequality  and \reff{eq:supMi-bd}
for the second inequality. 
Finally, from \eqref{eq:T2-0}, \eqref{eq:T2-1} and \eqref{eq:T2-2} (with
$\lim_{n\rightarrow \infty} \sum_{\ell = \p+1}^{\infty} (2\alpha)^{-\ell}
= 0$) , we get the result of the lemma.
\end{proof}

\section{Moments formula for BMC}
\label{sec:moment-BMC}
Let $X=(X_i,  i\in \T)$ be a  BMC on $(S, \cs)$  with probability kernel
$\cp$.          Recall           that          $|\G_n|=2^n$          and
$M_{\G_n}(f)=\sum_{i\in \G_n} f(X_i)$. We also recall that
$2\cq(x,A)=\cp(x, A\times S) + \cp(x, S\times A)$ for $A\in \cs$. 
We use the convention that $\sum_\emptyset=0$.
We recall the following well known and easy to establish  many-to-one
formulas  for BMC.

\begin{lem}
   \label{lem:Qi}
Let $f,g\in \cb(S)$, $x\in S$ and $n\geq m\geq 0$. Assuming that all the
quantities below are well defined,  we have:
\begin{align}
   \label{eq:Q1}
\E_x\left[M_{\G_n}(f)\right]
&=|\G_n|\, \cq^n f(x)= 2^n\, \cq^n f(x) ,\\
   \label{eq:Q2}
\E_x\left[M_{\G_n}(f)^2\right]
&=2^n\, \cq^n (f^2) (x) + 
 \sum_{k=0}^{n-1} 2^{n+k}\,   \cq^{n-k-1}\left( \cp
  \left(\cq^{k}f\otimes 
    \cq^k f \right)\right) (x),\\
   \label{eq:Q2-bis}
\E_x\left[M_{\G_n}(f)M_{\G_m}(g)\right]
&=2^{n} \cq^{m} \left(g \cq^{n-m} f\right)(x)\\
  \nonumber &\hspace{2cm} + \sum_{k=0}^{m-1} 2^{n+k}\, \cq^{m-k-1}
  \left(\cp\left(\cq^k g \sot \cq^{n-m+k} f\right) \right)(x). 
\end{align}
\end{lem}

We  also give  an upper  bound of  $\E_x\left[M_{\G_{n}}(f) ^4\right]$,
which is  a direct consequence  of the arguments  given in the  proof of
Theorem 2.1 in \cite{BDG14}.  Recall that $g\otimes^2=g\otimes g$. 
\begin{lem}
   \label{lem:M4}
   There exists a  finite constant $C$ such that for  all $f\in \cb(S)$,
   $n\in \N$ and  $\nu$ a probability measure on $S$,  assuming that all
   the quantities below are well defined, there exist functions $\psi_{j,
     n}$ for $1\leq j\leq 9$ such that:
\[
\E_\nu\left[M_{\G_n}(f)^4\right]= \sum_{j=1}^9 \langle \nu, \psi_{j, n}
\rangle,
\]
and,   with
$h_{k}= \cq^{k - 1} (f) $ and (notice that either $|\psi_j|$ or
$|\langle \nu, 
\psi_j \rangle|$ is bounded), writing  $\nu   g=\langle
\nu  ,   g   \rangle$:
\begin{align*}
 | \psi_{1, n}|
&\leq C \,2^n \cq^n(f^4),\\
 | \nu \psi_{2, n}|
&\leq C\,  2^{2n}\,
 \sum_{k=0}^{n-1} 2^{-k} |\nu \cq^k
 \cp  \left( \cq^{n-k - 1}( f^3) \sot h_{n- k} \right)|,\\
|\psi_{3, n}|
&\leq C 2^{2n} \sum_{k=0}^{n-1} 2^{-k}\,  \cq^k
  \cp \left( \cq^{n-k - 1} (f^2) \otimes^2
  \right),\\
|\psi_{4, n}|
&\leq C \, 2^{4n} \, 
\cp \left( |\cp(h_{n-1}\otimes^2)\otimes^2|\right), \\
|\psi_{5, n}|
&\leq C\,  2^{4n} \, 
 \sum_{k=2}^{n-1} \sum_{r=0}^{k -1}  2^{-2k-r }  \cq^r
  \cp \left( \cq^{k -r- 1} |\cp (h_{n- k} \otimes^2)|
  \otimes^2 \right),\\
|\psi_{6, n}|
&\leq C\, 2^{3n} \, 
  \sum_{k=1}^{n-1} \sum_{r=0}^{k -1} 
  2^{-k-r }  \cq^r| \cp \left(
\cq^{k -r-1}\cp \left( h_{n-k} \otimes^2 \right)\sot
\cq^{n-r-1}(f^2) \right)|,\\
|\nu \psi_{7, n}|
&\leq  C\, 2^{3n} \,   \sum_{k=1}^{n-1} \sum_{r=0}^{k -1} 
  2^{-k-r } |\nu \cq^r \cp \left(
\cq^{k -r-1}\cp \left( h_{n-k} \sot  \cq^{n-k -1} (f^2) \right)\sot
h_{n-r} \right)|,\\
|\psi_{8, n}|
&\leq C\,  2^{4n} \, 
  \sum_{k=2}^{n-1} \sum_{r=1}^{k -1} \sum_{j=0}^{r-1}
  2^{-k-r-j } \cq^j \cp \left(|
\cq^{r-j-1}\cp \left( h_{n-r} \otimes^2 \right)|\sot
|\cq^{k-j-1}\cp \left( h_{n-k} \otimes^2 \right)|
\right),\\
|\psi_{9, n}|
& \leq C\,  2^{4n} \, 
  \sum_{k=2}^{n-1} \sum_{r=1}^{k -1} \sum_{j=0}^{r-1}
  2^{-k-r-j } \cq^j |\cp \left(
\cq^{r-j-1}|\cp \left( h_{n-r} \sot \cq^{k -r -1}
  \cp\left(h_{n-k}\otimes^2 \right)\right) \sot h_{n-j} 
\right)|.
\end{align*}
\end{lem}

\bibliographystyle{abbrv}
\bibliography{biblio}

\end{document}